\documentclass[12pt]{amsart}
\usepackage{amssymb,amsthm}

\numberwithin{equation}{section}
 \newtheorem{thm}[equation]{Theorem}
\newtheorem{theorem}[equation]{Theorem}
 \newtheorem{defn}[equation]{Definition}
 \newtheorem{prop}[equation]{Proposition}

\newtheorem{cor}[equation]{Corollary}
 \newtheorem{lemma}[equation]{Lemma}

 \theoremstyle{definition}

 \newtheorem{remark}[equation]{Remark}

\DeclareMathOperator{\Ext}{Ext}
\DeclareMathOperator{\Hom}{Hom}
\DeclareMathOperator{\codim}{codim}
\DeclareMathOperator{\Ker}{Ker}
\DeclareMathOperator{\Ima}{Im}

\newcommand{\DOT}{\setlength{\unitlength}{1pt}\begin{picture}(2.5,2)
                  (1,1)\put(2,3.5){\circle*{2}}\end{picture}}

\newcommand{\HH}{{\rm HH}}
\newcommand{\Wedge}{\textstyle\bigwedge}
\newcommand{\CC}{\mathbb{C}}

\newcommand{\ZZ}{\mathbb{Z}}
\newcommand{\Sym}{\mathfrak S}
\newcommand{\ot}{\otimes} 

\newcommand{\A}{\mathcal{A}}

\title[Hochschild cohomology and graded Hecke algebras]
{Hochschild cohomology and graded Hecke algebras 
}

\date{July 18, 2006}
\author{Anne V.\ Shepler}
\address{Department of Mathematics, University of North Texas,
Denton, Texas 76203, USA}
\email{ashepler@unt.edu}
\thanks{The first author was partially supported by NSF grant
\#DMS-0402819. }
\author{Sarah Witherspoon}
\address{Department of Mathematics\\Texas A\&M University\\
College Station, Texas 77843, USA}\email{sjw@math.tamu.edu}
\thanks{The second author was partially supported by NSF grant
\#DMS-0443476 and the Alexander von Humboldt Foundation.}

\dedicatory{We dedicate this article to Sergey Yuzvinsky on
the occasion of his 70-th birthday.}

\begin{document}
\maketitle

\begin{abstract}
We develop and collect techniques for determining Hochschild
cohomology of skew group algebras $S(V)\# G$
and apply our results to graded Hecke algebras.
We discuss the explicit computation of certain types of invariants
under centralizer subgroups, focusing on the infinite family of complex
reflection groups $G(r,p,n)$ to illustrate our ideas. Resulting
formulas for Hochschild two-cocycles give information about deformations
of $S(V)\# G$ and, in particular, about graded Hecke
algebras.  We expand the definition
of a graded Hecke algebra to allow a nonfaithful action of $G$ on $V$,
and we show that there exist
nontrivial graded Hecke algebras for
$G(r,1,n)$, in contrast to the case of the natural reflection representation.
We prove that one of these graded Hecke algebras is equivalent to an
algebra that has appeared before in a different form.
\end{abstract}

%%%%%%%%%%%%%%%%%%%%%%%%%%%%%%%%%%%%%%%%%%%%%%%%%%%%%%%%%%%%%%%%%%%%%
\section{Introduction}

Lusztig~\cite{Lusztig3} showed that the graded version of an affine Hecke algebra
essentially retains its  representations,
which in turn determine certain representations of the corresponding group of Lie type 
\cite{KazhdanLusztig}.
Lusztig's graded Hecke algebra for a Coxeter group is a deformation of the skew group algebra
generated by the group and the vectors upon which the group acts
(in its natural reflection representation).  
Drinfeld~\cite{Drinfeld} had considered a different type of
deformation, a degenerate affine Hecke algebra: He allowed
{\em any} finite subgroup of ${\rm GL}(V)$ and 
used skew-symmetric forms on $V$ to deform the skew group algebra.
Ram and the first author \cite{RamShepler} showed that
Lusztig's graded Hecke algebra was a special case of Drinfeld's construction.
Drinfeld's algebra thus serves as a graded Hecke algebra for arbitrary
finite subgroups of ${\rm GL}(V)$.  One complication arises from this viewpoint
when working with a wider class of groups:
Drinfeld's construction for a given finite subgroup of ${\rm GL}(V)$ may only 
yield the trivial graded Hecke algebra.  This article is partly motivated
by a desire to understand and resolve this problem.  
We use Hochschild cohomology as a tool to explore those deformations of skew
group algebras that satisfy an expanded definition of graded Hecke algebra.
To illustrate our ideas (and as an application of results), we focus on the infinite family of complex reflection groups.

Drinfeld's construction yields only the trivial algebra for the complex reflection groups $G(r,1,n)=(\ZZ/r\ZZ)\wr\Sym_n$ when $r\geq 3$ and $n\geq 4$ 
(see \cite{RamShepler}).
Ariki and Koike \cite{ArikiKoike} defined a Hecke algebra (the Ariki-Koike algebra) for
these complex reflection groups by generators and relations.
Around the same time, Brou\'e and Malle \cite{BroueMalle} defined a cyclotomic Hecke algebra
for $G(r,1,n)$ (coinciding with the Ariki-Koike algebra) and
Brou\'e, Malle, and Rouquier~\cite{BroueMalleRouquier}
gave a topological description using braid groups
(see also \cite{BroueMalleMichel}).
Although the Hecke algebra itself enjoys fruitful study (see the survey paper by Mathas \cite{Mathas}),
a satisfactory formulation of an affine Hecke algebra for the groups $G(r,1,n)$
seems to be lacking.
One asks when a graded Hecke algebra may be a reasonable substitute.
In this article, we use Hochschild cohomology and algebraic deformation theory
to address a natural question: 
When and what kinds of graded Hecke algebras may be defined for $G(r,1,n)$?

Although Drinfeld's construction only yields trivial graded Hecke 
algebras for $G(r,1,n)$ in general, Ram and the first author \cite{RamShepler} 
constructed a novel graded Hecke algebra for these groups using ad hoc methods.
Dez\'el\'ee \cite{Dezelee} independently discovered the same algebra as a
subalgebra of a rational Cherednik algebra in the case $r=2$, that is, for the
Weyl group ${WB}_n = G(2,1,n)$, and defined this algebra for $WD_n=G(2,2,n)$.
Dez\'el\'ee \cite{Dezelee2} further realized the graded Hecke algebras of Ram
and the first author, for general $r$, 
as subalgebras of rational Cherednik algebras (equivalently,
graded Hecke algebras associated to
$V\oplus V^*$, where $V$ is the reflection representation).
These rational Cherednik algebras are examples of symplectic reflection algebras,
defined by Etingof and Ginzburg \cite{EtingofGinzburg} who were interested in
related deformations of orbifolds.

In this paper, we expand the definition of graded Hecke algebra 
to allow nonfaithful representations of a finite group $G$ 
on a finite dimensional vector space $V$ (Definition \ref{S:GHA}). 
We construct a natural graded Hecke algebra for $G(r,1,n)$
using a nonfaithful representation, and we
show that this construction coincides with the algebra
of Ram and the first author \cite{RamShepler} (Theorem~\ref{S:iso}).
We further give the generic form of a graded Hecke algebra associated 
to $G(r,1,n)$, depending on many parameters, in (\ref{S:general}).
A related paper of Chmutova \cite{Chmutova}
deals with existence of these more general
types of graded Hecke algebras defined for any symplectic representation
of $G$.

A graded Hecke algebra for a group $G$, with a given representation on
a finite dimensional vector space $V$, is a
particular type of deformation of the skew group algebra $S(V)\# G$.
Thus every graded Hecke algebra defines an
element in degree 2 Hochschild cohomology $\HH^2(S(V)\# G)$.
We characterize all possible graded Hecke algebras, both in relation to
Hochschild cohomology (Theorem \ref{S:paramspace}) and in relation
to defining skew-symmetric forms (Corollary \ref{S:ghacor}).
Corollary \ref{S:ghacor} generalizes a result of Ram and the first author
\cite{RamShepler}, first formulated by Drinfeld \cite{Drinfeld}
for Coxeter groups.

There may well be other compelling  deformations of $S(V)\# G$
(arising from some of the many elements of degree 2
Hochschild cohomology) that are not graded Hecke algebras.
Some small examples of such deformations were given in 
\cite{CaldararuGiaquintoWitherspoon,Witherspoon1}.
With these potential deformations in mind, as well as the immediate
applications to graded Hecke algebras in this article,
we compute the relevant cohomology in two cases:
\begin{itemize}
\item[(1)] $G$ is the monomial reflection group $G(r,p,n)$
acting via its natural reflection representation on $V=\CC^n$
(Theorem \ref{main-faithful-theorem}), and
\item[(2)] $G=G(r,p,n)$ acting on $V=\CC^n$ via  the permutation
representation of its quotient symmetric group $\Sym_n$ 
(Theorem \ref{main-nonfaithful-theorem}).
\end{itemize}

The graded vector space structure of Hochschild cohomology 
$\HH^{\DOT}(S(V)\# G)$ can be described
in terms of invariants of
centralizer subgroups of $G$ (see \cite{Farinati,GinzburgKaledin}).
We recall this structure in Section \ref{terminologyandnotation} and
reformulate it in terms of semi-invariants.
In Sections 5--7, we explicitly compute these semi-invariants in degree 2 using
techniques of classical invariant theory and the theory of hyperplane arrangements;
the full Hochschild cohomology space $\HH^{\DOT}(S(V)\# G)$ could be
computed similarly  using our approach.
Although we are primarily interested in applications to graded Hecke
algebras and other deformations, our explicit computation of Hochschild 
cohomology should also be of interest to homological algebraists and to geometers
studying cohomology for associated orbifolds.

Our results in case (1) above show that the lack of nontrivial graded Hecke
algebras (when $r\geq 3$, $n\geq 4$, see \cite{RamShepler}) 
stems from a lack of
certain types of elements in $\HH^2(S(V)\# G)$. 
This conclusion 
motivates our computation in case (2), 
leading to the aforementioned nontrivial graded Hecke algebras 
under a nonfaithful representation of $G$.

The authors thank A.\ Ram for questions and suggestions that led to
this project.

%%%%%%%%%%%%%%%%%%%%%%%%%%%%%%%%%%%%%%%%%%%%%%%%%%%%%%%%%%%%%%%%%%%%

%%%%%%%%%%%%%%%%%%%%%%%%%%%%%%%%%%%%%%%%%%%%%%%%%%%%%%%%%%%%5
\section{Terminology and notation}

We work over the complex numbers $\CC$, letting $\otimes
=\otimes_{\CC}$ unless otherwise specified.
Let $G$ be a finite group and let $V$ be the vector space $\CC^n$ with
a representation of $G$.
Let $V^*$ denote the contragredient (or dual) representation.
We will often need a choice of basis $v_1,\ldots,v_n$ of $V$ and
dual basis $x_1,\ldots,x_n$ of $V^*$.

Let $R$ be an associative $\CC$-algebra with an action of $G$ by
automorphisms.
In this paper, $R$ will always be either the symmetric algebra
$S(V)$ or the tensor algebra $T(V)$.
The {\bf skew group algebra} $R\# G$ is the vector space
$R\ot \CC G$ with multiplication
$$
  (r\ot g)(s\ot h)=r\cdot g(s) \ot gh
$$
for all $r,s\in R$ and $g,h\in G$.
We abbreviate $r\ot g$ by $r\overline{g}$ ($r\in R$, $g\in G$)
and $r\ot 1$ or $1\ot g$ simply by $r$ or $\overline{g}$, respectively.
Note that an element $g\in G$ acts on $R$ by 
conjugation by $\overline{g}$: 
$ \ \overline{g} r (\overline{g})^{-1} =g(r)\overline{g} (\overline{g})
^{-1} =g(r)$ for all $r\in R$.

For our cohomological computations, we fix some
notation involving the action of $G$ on $V$.
Let $V^G=\{v \in V: g(v)=v \text{ for all } g\in G\}$, the set of 
$G$-invariants in $V$. More generally,
let $\chi :G\rightarrow \CC$ be any linear character (i.e., any
group homomorphism from $G$ to $\CC ^{\times}$) and let
$V^\chi=\{v\in V:g(v)=\chi(g)v \text{ for all } g\in G\}$, the set of 
{\bf $\chi$-invariants} in $V$ (also called {\bf semi-invariants} with 
respect to $\chi$).

For any $g\in G$, let $Z(g)=\{h\in G: gh=hg\}$, 
the centralizer of $g$ in $G$,
and let $V^g = \{v\in V: g(v)=v\}$, the $g$-invariant subspace of $V$.
Since $G$ is finite, we may assume $G$ acts by isometries on $V$ (i.e.\
$G$ preserves a Hermitian form).
If $h\in Z(g)$, then $h$ preserves both $V^g$ and its orthogonal 
complement $(V^g)^\perp$ and we  define
$$ 
 h^\perp = h|_{(V^g)^\perp}
$$
(i.e., $h^{\perp}$ is the linear transformation by which $h$ acts on
the vector space $(V^g)^\perp$).
We are particularly interested in semi-invariants with respect to
the {\bf Hochschild character}
$\chi_g:Z(g)\rightarrow \CC$ 
{\bf of} $g$ defined by
\begin{equation}\label{S:hochchar}
  \chi_g(h)= \det(h^{\perp})
\end{equation}
for all $h\in Z(g)$.  
When no confusion will arise, we simply 
write $\chi$ instead of $\chi_g$.
Note that $\chi_g$ is a linear character.
Let $1$ denote the identity element of $G$ 
and let $I_n$ (or just $I$) denote the $n \times n$ identity matrix. 

Our results involve the following complex reflection groups
parametrized by positive integers $r,p,n$:
The {\bf full monomial group} $G(r,1,n) < \text{GL}_n(\CC)$ 
is the group of all monomial $n\times n$ matrices 
whose nonzero entries are $r$-th roots of unity.
The full monomial group is a finite complex reflection group 
acting on $V=\CC^n$ and is isomorphic
to the wreath product of
the cyclic group of order $r$ and the symmetric group:
$$
 G(r,1,n) = N\cdot G(1,1,n) 
      \cong  (\ZZ/r\ZZ)^n \rtimes \Sym_n
      = \ZZ/r\ZZ \wr \Sym_n , 
$$
where $N$ is the normal subgroup of diagonal matrices
in $G(r,1,n)$ and $G(1,1,n)\leq G(r,1,n)$ is
isomorphic to the symmetric group $\Sym_n$.

Let $\xi$ in $\CC$ be a fixed primitive $r$-th root of unity.
Let $\xi_i$ be the diagonal 
matrix $\text{diag}\{1, \ldots, 1, \xi, 1, \ldots, 1\}$,
with $i$-th entry $\xi$.  
We denote each element in $G(r,1,n)$ 
by $g=\xi_1^{a_1} \ldots\xi_n^{a_n} \sigma$,
where $\sigma \in \Sym_n$ and $0 \leq a_i < r$.
In the special case that $r$ is even (so $\xi^{r/2} = -1$),
we write
 $(1, -2):=\xi_2^{r/2}(1,\, 2)$ for the transformation
$v_2 \mapsto v_1$, $v_1 \mapsto -v_2$.

For any $p$ dividing $r$, we consider a natural subgroup of $G(r,1,n)$
also acting as a reflection group on $V$.
Let $\psi(M)$
be the product of the nonzero entries of any matrix $M$.  Then
$$
 G(r,p,n)=\{M\in G(r,1,n):\ \psi(M)^{r/p}=1 \},
$$
that is,
$M$ in $G(r,1,n)$ lies in $G(r,p,n)$ when the sum of the exponents
of $\xi$ in the nonzero entries of $M$ is $0$ modulo $p$.  
Several real reflection groups (Coxeter groups) are special cases 
of the groups $G(r,p,n)$:
\begin{itemize}
\item $G(1,1,n)$ is the symmetric group $\Sym_n$,
\item $G(2,1,n)$ is the Weyl group ${WB}_n$ of type $B_n$,
\item $G(2,2,n)$ is the Weyl group ${WD}_n$ of type $D_n$, and 
\item $G(r,r,2)$ is the dihedral group $I_2(r)$ of order $2r$.
\end{itemize}

We briefly record some facts on centralizers and conjugacy classes in $G(r,1,n)$.
Note that $Z_{G(r,p,n)}(g)=Z_{G(r,1,n)}(g)\cap G(r,p,n)$ for any $g$ in $G(r,p,n)$. The (generalized) cycle structure of an element $g$ in $G(r,1,n)$ determines its conjugacy class and the size of its centralizer $Z(g)$ in the following way.
Every element $g$ in $G(r,1,n)$ can be written as the product of disjoint cycles of the form 
$\xi_i^{a_i} \xi_{i+1}^{a_{i+1}}\cdots\xi_{m}^{a_{m}}(i, i+1,  \ldots ,m)$.
We call such an element an $(a, k)$-cycle, where $k=m-i+1$ and
$a=a_i+a_{i+1}+\cdots +a_{m} \mod r$.
Two elements in $G(r,1,n)$ are conjugate if and only if they have the same number of $(a,k)$-cycles for each $a$ and $k$ (where $0\leq a < r$ and $1\leq k \leq n$).
Let $m_{a,k}$ be the number of $(a,k)$-cycles for $g$;
the order of the centralizer of $g$ in $G(r,1,n)$ is
\begin{equation}
\label{centralizer}
|Z(g)| = \prod_{a,k} m_{a,k}! \cdot k^{m_{a,k}} \cdot r^{m_{a,k}}.
\end{equation}
See
\cite[Equation (2.5)]{RamShepler} and \cite[Appendix B, Section 3]{Macdonald}.

%%%%%%%%%%%%%%%%%%%%%%%%%%%%%%%%%%%%%%%%%%%%%%%%%%%%%%%%%%%%%%%%%%%%%%%
%%%%%%%%%%%%%%%%%%%%%%%%%%%%%%%%%%%%%%%%%%%%%%%%%%%%%%%%%%%%%%%%%%%%%%%
%%%%%%%%%%%%%%%%%%%%%%%%%%%%%%%%%%%%%%%%%%%%%%%%%%%%%%%%%%%%%%%%%%%%%%%
%%%%%%%%%%%%%%%%%%%%%%%%%%%%%%%%%%%%%%%%%%%%%%%%%%%%%%%%%%%%%%%%%%%%%%%
%%%%%%%%%%%%%%%%%%%%%%%%%%%%%%%%%%%%%%%%%%%%%%%%%%%%%%%%%%%%%%%%%%%%%%%
\section
{Hochschild cohomology of $S(V)\# G$}
\label{terminologyandnotation}

The {\bf Hochschild cohomology} of a $\CC$-algebra $R$ is the graded
vector space
$\HH^{\DOT}(R)=\Ext^{\DOT}_{R^e}(R,R)$, where $R^e=R\ot R^{op}$ acts
on $R$ by left and right multiplication.
More generally,
if $M$ is an $R$-bimodule (equivalently a left $R^e$-module),
$\HH^{\DOT}(R,M)=\Ext^{\DOT}_{R^e}(R,M)$ so that 
$\HH^{\DOT}(R)=\HH^{\DOT}(R,R)$.
For more details, see \cite{Weibel}.

Fix a finite group $G$ and a representation $V=\CC^n$ of $G$. 
The Hochschild cohomology of $S(V)\# G$ is given in 
\cite{Farinati,GinzburgKaledin} when $G$ acts faithfully on $V$.
In this section, we reformulate this result
to aid our explicit computations and our determination of graded Hecke algebras.
We state all results more generally for representations which are not
necessarily faithful.
No new techniques are needed for the generalization, and so we merely
sketch those arguments that appear elsewhere for the
faithful case.

Let $\mathcal C$ be  a set of representatives of the conjugacy classes
of $G$. In the following theorem, 
a negative exterior power is interpreted as $0$.
Recall the Hochschild characters $\chi_g$ defined in (\ref{S:hochchar}).

\begin{thm}\label{S:hhchi}
There is an isomorphism of graded vector spaces
\begin{equation*}
  \HH^{\DOT}(S(V) \# G) \cong
  \bigoplus_{
   \substack{g \in {\mathcal C}}} 
  \left(S(V^g) \otimes \Wedge^{\DOT -\codim V^g} ((V^g)^*)\right)^{\chi_g}.
\end{equation*}
\end{thm}

For any $g\in{\mathcal {C}}$, let
$\HH^2(g)$ be the $g$-component of $\HH^2 (S(V)\# G)$ in 
the theorem:
\begin{equation}\label{S:hh2g}
  \HH^2(g)= 
  \left( S(V^g) \otimes \Wedge^{2 - \codim V^g} ((V^g)^*) \right)^{\chi_g}.
\end{equation}
We call $\HH^2(g)$ the set of {\bf Hochschild semi-invariants of} $g$.
Note $\HH^2(g)=0$ if $\codim V^g >2$, since the exterior
power $2\! -\!\codim V^g$ is negative in this case.

\begin{proof}[Sketch of a proof of Theorem \ref{S:hhchi}]
We give only an outline as the bulk of the calculation
follows \cite{Farinati,Witherspoon1}.
A result of \c{S}tefan on Hopf Galois extensions \cite[Cor.\ 3.4]{Stefan}
implies that there is an action of $G$ on $\HH^{\DOT}(S(V),S(V)\# G)$
for which 
$$
  \HH^{\DOT}(S(V)\# G)\cong \HH^{\DOT}(S(V), S(V)\# G)^G.
$$
Now $S(V)\# G\cong \oplus_{g\in G} S(V)\overline{g}$ as an
$S(V)$-bimodule.
The above isomorphism of graded vector spaces may thus be rewritten as
\begin{equation}\label{S:reynolds}
  \HH^{\DOT}(S(V)\# G)\cong \left(\bigoplus_{g\in G}\HH^{\DOT}(S(V),
      S(V) \overline{g}))\right)^G.
\end{equation}
The action of $G$ permutes the direct summands by conjugation on their
indices $g\in G$.
We may therefore rewrite the vector space structure of these $G$-invariants as
$\oplus_{g\in {\mathcal C}} \HH^{\DOT}(S(V), S(V)\overline{g} )^{Z(g)}$
(one applies a transfer (trace) operator $\sum_{h\in G/Z(g)} h$ to each 
$Z(g)$-invariant to obtain a $G$-invariant).
To compute the summands $\HH^{\DOT}(S(V),S(V)\overline{g})^{Z(g)}$,
one may first use a Koszul complex (see \cite[\S4.5]{Weibel} for details
of Koszul complexes) to find $\HH^{\DOT}
(S(V),S(V)\overline{g})$ and then determine the $Z(g)$-invariants.
The use of a Koszul complex yields
\begin{equation}\label{S:hhsg}
  \HH^{\DOT}(S(V)\# G)\cong \bigoplus_{g\in {\mathcal C}}
  \left(S(V^g)\overline{g}\ot
  \Wedge^{\DOT - \codim V^g}((V^g)^*)\ot \Wedge^{\codim V^g}(((V^g)
  ^{\perp})^*)\right)^{Z(g)},
\end{equation}
where $Z(g)$ acts diagonally on the given tensor product of 
three representations.
Its action on the factor $S(V^g)\overline{g}$ is by conjugation
by $\overline{h}$ ($h\in Z(g)$).
The actions on the given dual vector spaces are the contragredient actions.
The factor $\Wedge^{\codim V^g}(((V^g)^{\perp})^*)$ is one-dimensional
since $\codim V^g = \dim ((V^g)^{\perp})^*$.
We also denote this factor by $\det (((V^g)^{\perp})^*)$; it is
included in the above expression (\ref{S:hhsg})
only since it carries a potentially
nontrivial $Z(g)$-action.

We find the contribution to Hochschild cohomology
from each $g\in G$ representing a conjugacy class by 
determining the $Z(g)$-invariants of
\begin{equation}\label{gcomponent}
  S(V^g)\overline{g}\ot\Wedge^{\DOT -\codim V^g} ((V^g)^*)
  \otimes \det(((V^g)^{\perp})^*).
\end{equation}
Let $x_1,\ldots,x_m$ be a basis of $((V^g)^{\perp})^*$. 
Let $s\in S(V^g)\ot \Wedge^{\DOT - \codim V^g}((V^g)^*)$.
Then the element $\overline{g}s\ot x_1\wedge\cdots\wedge x_m$
of $S(V^g)\overline{g}\ot\Wedge^{\DOT -\codim V^g} ((V^g)^*)
\otimes \det(((V^g)^{\perp})^*)$ is $Z(g)$-invariant if and only if
$$  \overline{g}s\ot x_1\wedge\cdots\wedge x_m 
 = \det((h^{-1})^{\perp})\, \overline{g}\, h(s)\ot x_1\wedge
  \cdots \wedge x_m
$$
for all $h\in Z(g)$.
Equivalently,
$
  h(s) = \det(h^{\perp}) s =\chi_g(h)s
$
for all $h\in Z(g)$.
\end{proof}

The next lemma gives in particular a necessary condition for $\HH^2(g)$
to be nonzero (cf.\ \cite[Example 3.10]{Farinati}).

\begin{lemma}
\label{forceszero}
Let   $g\in G$, and suppose $(S(V^g)\otimes \Wedge^{i -\codim V^g}
((V^g)^*))^{\chi_g}\neq 0$ for some $i\geq 0$. 
If $h\in Z(g)$ and $h|_{V^g}$ is the identity map, then $\det(h)=1$.
In particular, $\det (g) =1$.
\end{lemma}

\begin{proof}
Let $s$ be a nonzero element of 
$(S(V^g)\otimes \Wedge^{i -\codim V^g}((V^g)^*))^{\chi_g}$.
If $h\in Z(g)$ and $h|_{V^g}$ is the identity,
then $h(s)=s$, and so
$\det(h) s= \det(h^\perp) s= \chi_g(h) s= h(s) =  s$.
\end{proof}

We may now rewrite Theorem \ref{S:hhchi} 
for degree 2 cohomology using (\ref{S:hh2g}) and the above lemma.
The condition $\det(g)=1$ from the lemma implies $\codim V^g \neq 1$
if $\HH^2(g)\neq 0$, and the
appearance of the exterior power $2\! -\!\codim V^g$ in (\ref{S:hh2g}) implies
that $\codim V^g \leq 2$.
Thus we have
\begin{equation}
\label{S:det1}
 \HH^{2}(S(V) \# G) \cong
  \bigoplus_{\substack{g \in {\mathcal C}}} \HH^2(g) \cong
  \bigoplus_{\substack{g \in {\mathcal C}\\\rule[0ex]{0ex}{1.5ex}%strut
  \det(g)=1\\ \rule[0ex]{0ex}{1.5ex}%strut
   \codim V^g\in\{0,2\}}} \HH^2(g).
\end{equation}

\begin{remark}
We briefly mention a generalization that we do not pursue in this paper.
Let $\alpha:G\times G\rightarrow \CC^{\times}$ be a two-cocycle, that is
$\alpha(g,h)\alpha(gh,\ell)=\alpha(h,\ell)\alpha(g,h\ell)$ for all
$g,h,\ell\in G$.
The {\bf crossed product algebra} $S(V)\#_{\alpha}G$ is the vector space 
$S(V)\ot\CC G$ with multiplication $(r\overline{g})
(s\overline{h}) = \alpha(g,h) r\cdot g(s)\overline{gh}$ for all
$r,s\in R$, $g,h\in G$.
The derivation of the more general
Hochschild cohomology $\HH^{\DOT}(S(V)\#_{\alpha}G)$ is 
similar to that of $\HH^{\DOT}(S(V)\# G)$.
(See \cite[Cor.\ 6.5]{Witherspoon1}, which remains valid for
nonfaithful actions of $G$ on $V$.)
One obtains a sum as in Theorem \ref{S:hhchi}, namely
\begin{equation*}
  \HH^{\DOT}(S(V) \#_{\alpha} G) \cong
  \bigoplus_{
   \substack{g \in {\mathcal C}}} 
  \left(S(V^g) \otimes \Wedge^{\DOT -\codim V^g} ((V^g)^*)\right)^{\chi
    ^{\alpha}_g}.
\end{equation*}
However, the linear character $\chi^{\alpha}_g$ on $Z(g)$ now depends on $\alpha$:
$$
  \chi^{\alpha}_g(h)=\det (h^{\perp})\alpha(g,h)(\alpha(h,g))^{-1}
$$
for all $h\in Z(g)$.
These potentially different characters $\chi^{\alpha}_g$ can
lead to different semi-invariants.
In this broader setting, deformations of $S(V)\#_{\alpha}G$ are
treated in 
\cite{CaldararuGiaquintoWitherspoon,Chmutova,Witherspoon1,Witherspoon2}.
\end{remark}

%%%%%%%%%%%%%%%%%%%%%%%%%%%%%%%%%%%%%%%%%%%%%%%%%%%%%%%%%
%%%%%%%%%%%%%%%%%%%%%%%%%%%%%%%%%%%%%%%%%%%%%%%%%%%%%%%%%%%%%%%%%%%%%%%
%%%%%%%%%%%%%%%%%%%%%%%%%%%%%%%%%%%%%%%%%%%%%%%%%%%%%%%%%%%%%%%%%%%%%%%
%%%%%%%%%%%%%%%%%%%%%%%%%%%%%%%%%%%%%%%%%%%%%%%%%%%%%%%%%%%%%%%%%%%%%%%
%%%%%%%%%%%%%%%%%%%%%%%%%%%%%%%%%%%%%%%%%%%%%%%%%%%%%%%%%%%%%%%%%%%%%%%
%%%%%%%%%%%%%%%%%%%%%%%%%%%%%%%%%%%%%%%%%%%%%%%%%%%%%%%%%%%%%%%%%%%%%%%
\section{Invariant Theory of $G(r,p,n)$}
\label{invthyofG(r,p,n)}
 
The Hochschild cohomology of the skew group algebra 
$S(V)\# G(r,p,n)$ is derived from
the invariant theory of $G(r,p,n)$ and its centralizer subgroups.
In this section, we briefly record the invariant theory of $G(r,p,n)$
acting via its natural reflection representation on $V=\CC^n$
and on the exterior algebra of derivations on $V^*$.
We first recall some well-known facts on reflection groups.

Recall
 that a {\bf reflection} is an element of $\text{GL}(V)$ of finite order which
fixes a hyperplane (called the {\bf reflecting hyperplane}) in $V$ pointwise.
Let $G$ be a reflection group, i.e., a (finite) group $G\leq \text{GL}(V)$ generated by reflections.  The collection of reflecting hyperplanes $\A$ for $G$ is called the {\bf reflection arrangement}.  For each $H\in \A$,
let $l_H$ in $V^*$ be a linear form defining $H =\ker l_H$.
The reflection arrangement $\A$ is 
defined (up to a nonzero constant) by the polynomial
$Q = \prod_{H\in\A} l_H$ in $S(V^*)$.
When $G$ is a reflection group acting on
$V^*$, we identify $(V^*)^*$ with $V$ and define the reflection arrangement
$\A$ by a polynomial $Q$ in $S(V)$.

%%%%%%%%%%%%%%%%%%%%%%%%%%%%%%%%%%%%%%%%%%
\begin{remark}
\label{directsumofreflectiongroups}
By the Shephard-Todd-Chevalley Theorem (see 
\cite[Chapter 6]{OrlikTerao}),
the invariant ring $S(V)^G$ is a polynomial algebra generated by $n$
algebraically independent invariants: $S(V)^G=\CC[f_1, \ldots, f_n]$,
for some polynomials $f_i$ called {\bf basic invariants}.
In fact, a set of homogeneous, algebraically independent, invariant polynomials $f_1,\ldots, f_n$ generate $S(V)^G$ if and only if the product of the degrees
$\deg f_i$ is $|G| $ (for example, see \cite[Theorem 3.7.5]{KemperDerksen}).
Hence, the invariants
of a direct sum of reflection representations is
the tensor product of the invariants:
If $G$ and $G'$ are reflection groups acting on complex vector spaces
$V$ and $V'$ with basic invariants
 $f_1, \ldots, f_n$ and
$f'_1, \ldots, f'_{n'}$, respectively, then
$$
 S(V \oplus V')^{G \oplus G'} 
  =\CC[f_1,\ldots, f_{n}, f'_1,\ldots, f'_{n'}]
  \cong S(V)^{G} \otimes S(V')^{G'}.
$$
\end{remark}
%%%%%%%%%%%%%%%%%%%%%%%%%%%%%%%%%%%%%%%%%%%%%%%%%%55

The following proposition (see \cite{ShephardTodd}, for example, or above remark)
describes the polynomials invariant under $G(r,p,n)$.
Recall that the $k$-th elementary
symmetric function in variables 
$z_1, \ldots, z_n$ is the polynomial
$f_k = \sum_{1\leq i_1 < \cdots < i_k \leq n} z_{i_1} \cdots z_{i_k}$.
%
%%%%%%%%%%%%%%%%%%%%%%%%%%%
\begin{prop}
\label{invariantsofG(r,p,n)} 
Let $G=G(r,p,n)$ act on $V$ with $\CC$-basis $v_1, \ldots, v_n$
via its natural reflection 
representation. Then
$$
 S(V)^{G(r,p,n)} = \CC[f_1, \ldots, f_n]
$$
where $f_i$ is the $i$-th elementary symmetric function in
variables $v_1^r, \ldots, v_n^r$ for $i<n$
and $f_n = (v_1 \cdots v_n)^{r/p}$.
\end{prop}
%%%%%%%%%%%%%%%%

We identify $S(V) \otimes  \Wedge^k (V^*)$
with the $S(V)$-module of $k$-derivations
(or $k$-vector fields) on $V^*$.
The following proposition (see \cite{OrlikTerao})
describes the invariant derivations on $V^*$.  The generators $\theta_i$
in the proposition are called {\bf basic derivations}. 

%%%%%%%%%%%%%%%%%%%%%%%%%%%%%%%%%%%%%%%
\begin{prop}
\label{Solomon'sThm}
Let $G$ be a reflection group acting on $V^*=\CC^n$ with reflection arrangement defined by $Q \in S(V)$.  Suppose $\theta_1,\ldots, \theta_n$
in $(S(V)\otimes V^*)^G$ are invariant derivations
whose coefficient matrix has determinant $Q$ up to a nonzero scalar.
Then
$$
\begin{aligned}
     \left(S(V) \otimes  \Wedge^k (V^*)\right)^{G}
   = 
     \bigoplus_{1\leq i_1 <\cdots < i_k \leq n}\ 
     S(V)^{G}\ (\theta_{i_1} \wedge\cdots\wedge \theta_{i_k}).
\end{aligned}
$$
\end{prop}
%%%%%%%%%%%%%%%%%%%%%%%%%%%%%%
\begin{proof}
Note that $ S(V)^{\det^{-1}} = S(V)^{G}\ Q$, 
where $\det^{-1}$ is the inverse of the determinant character
of $G$ acting on $V^*$ (see \cite[Example 6.40]{OrlikTerao}).
Solomon's Theorem \cite[Prop.~6.47]{OrlikTerao} 
then implies that the
invariant $k$-derivations are generated freely
over $S(V)^{G}$
by the nonzero wedge products of the $\theta_i$
(taken $k$ at a time).
\end{proof}
%%%%%%%%%%%%%%%%%%%%%%%%%%%%%%%%%%%%%%%%
%%%%%%%%%%%
%
We apply the above proposition to $G=G(r,p,n)$ (see \cite[Appendix B]{OrlikTerao}) .
%%%%%%%%%%%%%%%%%%%%%
\begin{prop}
\label{explicitderivations}
Let $G(r,p,n)$ act on $V=\CC^n$ with $\CC$-basis $v_1, \ldots, v_n$ 
and on $V^*$ contragadiently with dual basis $x_1, \ldots, x_n$.
Then
$$
\begin{aligned}
     \left(S(V) \otimes  \Wedge (V^*)\right)^{G(r,p,n)}
   = \bigoplus_{1\leq k\leq n}\
     \bigoplus_{1\leq i_1 <\cdots < i_k \leq n}\ 
     S(V)^{G(r,p,n)}\ (\theta_{i_1} \wedge\cdots\wedge \theta_{i_k}),
\end{aligned}
$$
where
$$
 \theta_j := \sum_{1\leq i \leq n} v_i^{(j-1)r+1} \otimes  x_i\\
$$
for $j<n$ and
$$
 \theta_n :=
   \begin{cases}
      \ \sum_{1\leq i\leq n}\ v_i^{(n-1)r+1} \otimes  x_i
      &\text{ if } p\neq r\\
      \ \sum_{1\leq i\leq n}\ (v_1\cdots v_{i-1} v_{i+1}\cdots v_n)^{r-1} \otimes  x_i
      &\text{ if } p= r 
   \end{cases}   \ \ \ .
$$
\end{prop}

%%%%%%%%%%%%%%%%%%%%%%%%%%%%%%%%%%%%%%%%%%%%%%%%%%%%%%%%%%%%%%%%%%%%%%%
%%%%%%%%%%%%%%%%%%%%%%%%%%%%%%%%%%%%%%%%%%%%%%%%%%%%%%%%%%%%%%%%%%%%%%%
%%%%%%%%%%%%%%%%%%%%%%%%%%%%%%%%%%%%%%%%%%%%%%%%%%%%%%%%%%%%%%%%%%%%%%%
%%%%%%%%%%%%%%%%%%%%%%%%%%%%%%%%%%%%%%%%%%%%%%%%%%%%%%%%%%%%%%%%%%%%%%%
%%%%%%%%%%%%%%%%%%%%%%%%%%%%%%%%%%%%%%%%%%%%%%%%%%%%%%%%%%%%%%%%%%%%%%%
\section{Hochschild cohomology in degree 2 for $G(r,p,n)$}
\label{faithfulHochschildcohomologyindeg2}

We give the Hochschild $2$-cohomology
$\HH^2\left( S(V)\# G \right)$ explicitly 
when the monomial group $G=G(r,p,n)$ acts 
(faithfully) via its natural reflection representation
on $V=\CC^n$.
The results rely on some computations postponed until Section 
\ref{semi-invariantsofG(r,p,n)actingfaithfully}.
We assume $n\geq 4$ throughout this section.  Results for $n=2$ and $n=3$ are slightly different because the pattern of relevant
conjugacy classes in $G$ differs for low values of $n$.  The techniques from this section apply to give explicit Hochshild cohomology in the case $n=2$ or $n=3$, but we omit these results for sake of brevity.  
In the following theorem, we use the notation $\HH^2(g)$ defined in 
(\ref{S:hh2g}).

%%%%%%%%%%%%%%%%%%%%%%%%%%%%%%%%%%%%%%%%%%%%%
\begin{theorem}\label{main-faithful-theorem}
Assume
$n\geq 4$.
Let $G=G(r,p,n)$ act on $V=\CC^n$ via the natural reflection
representation.
The Hochschild cohomology in degree $2$ for the skew group algebra
$S(V)\#G$ is 
$$
  \HH^2(S(V)\#G) \ \cong 
  \ \HH^2(1) \oplus \HH^2(1, 2,3)\oplus \HH^2(1, -2) \oplus 
     \bigoplus_{i=1,\ldots, \lfloor r/2 \rfloor}
    \HH^2(\xi_1^i\ \xi_2^{-i})
$$
where 
\begin{itemize}
\item
$\HH^2(1) =\left( S(V) \otimes \Wedge^{2} (V^*) \right)^{G(r,p,n)}$ 
is given explicitly in Proposition~\ref{explicitderivations},
\item
$\HH^2(1, 2, 3)$ is given explicitly in Proposition~\ref{HH^2(g_b)},
\item
$\HH^2(1, -2)=0$ unless $r=2p$,
in which case
$\HH^2(1, -2)$ is given explicitly in Proposition~\ref{HH^2(g_c)},
\item
$\HH^2(\xi_1^\ell\ \xi_2^{-\ell})=0$ unless $r>1$, $p=r$, and
$\ell\neq r/2$,
in which case
$\HH^2(\xi_1^\ell\ \xi_2^{-\ell})$ is given explicitly in Proposition~\ref{HH^2(g_d)3}.
\end{itemize}
\end{theorem}
%%%%%%%%%%%%%%%%%%%%%%%%%%%%%%%%%%%%%%%%%%%%%%%%%%%%%
\begin{proof}
We apply (\ref{S:det1}).
If $g\in G(r,p,n)$ with $\codim(V^g)=0$, then $g$ is the
identity element $1$ (as $G(r,p,n)$ acts faithfully on $V$).
Furthermore, the Hochschild character $\chi$
for $g=1$ should be interpreted as 
the trivial character of $Z(g)=G(r,p,n)$
since $(V^g)^\perp=\{0\}$.  Hence
$$
  \HH^2(1)= 
  \left( S(V^g) \otimes \Wedge^{2} ((V^g)^*) \right)^{G(r,p,n)},
$$
which is given explicitly in Proposition \ref{explicitderivations}.

If $g\in G(r,p,n)$ with $\codim(V^g)=2$,
then $g$ is conjugate to one of the following elements
(see \cite[Section 2B]{RamShepler}):
$$\begin{matrix}
g_b= \xi_1^{a}\xi_3^{-a}(1, 2, 3), \hfill
&\qquad &\hbox{$0\le a\le \gcd(p,3)-1$,} \hfill \cr
\cr
g_c=\xi_1^{a+\ell}\xi_2^{-a}(1, 2),\hfill 
&\qquad &\hbox{$\ell\neq 0$ (so $r\ne 1$),} \hfill \cr
\cr
g_d=\xi_1^{\ell_1}\xi_2^{\ell_2},\hfill 
&&\hbox{$\ell_1 \neq 0$, $\ell_2\ne 0$ (so $r\ne 1$)}, \hfill\cr
\cr
g_e=(1, 2)\xi_3^\ell,\hfill && \ell \neq 0, \hfill\cr 
\cr
g_f=(1, 2)\xi_3^a\xi_4^{-a}(3, 4). \cr
\end{matrix}
$$

We can further conjugate these elements into simpler forms under our
assumption $n\geq 4$.
Each element of the form
$g_b= \xi_1^{a} \xi_3^{-a}\ (1, 2, 3)$ is conjugate
to $(1, 2, 3)$ 
(via $(1, 2, 3)\xi_3^a\xi_4^{-a}$) in $G(r,p,n)$.
Each element of the form 
$g_c= \xi_1^{a+\ell} \xi_2^{-a}\ (1, 2)$ is conjugate
to $\xi_2^{\ell}(1, 2)$ (via $\xi_1^a(1, 2)\xi_3^{-a}$) 
in $G(r,p,n)$.
We assume $1=\det g_c= -\xi^\ell$, else $\HH^2(g_c)=0$ by
Lemma \ref{forceszero}. 
Thus, if $\HH(g_c)$ is nonzero, $r$ must be even with $\xi^\ell=-1$
and we may assume $g_c=\xi_2^\ell (1, 2)=(1, -2)$.
The sets of Hochschild semi-invariants $\HH^2(1, 2, 3)$ and
$\HH^2(1, -2)$ are given in Propositions~\ref{HH^2(g_b)} and \ref{HH^2(g_c)}.

By Proposition~\ref{HH^2(g_d)1}, 
the set of Hochschild semi-invariants $\HH^2(g_d)$
is zero unless $G=G(r,r,n)$.  The element $g_d$ lies in $G(r,r,n)$
only when $\ell_2\equiv -\ell_1 \mod r$.
In the special case that $\ell_1= \ell_2= r/2 $, $\HH^2(g_d)$ is zero by
Proposition \ref{HH^2(g_d)2}.  Otherwise,
$\HH^2(g_d)$ is given in Proposition~\ref{HH^2(g_d)3}.

The sets of Hochschild semi-invariants $\HH^2(g_e)$
and $\HH^2(g_f)$ are both zero by 
Propositions~\ref{HH^2(g_e)} and \ref{HH^2(g_f)}.
\end{proof}
%%%%%%%%%%%%%%%%%%%%%%

The Hochschild cohomology in degree 2 simplifies when
$G$ is the symmetric group $\Sym_n=G(1,1,n)$ and 
when $G$ is the Weyl group ${WB}_n=G(2,1,n)$.
(See Proposition~\ref{invariantsofG(r,p,n)} for the explicit rings of invariants
$S(V)^{\Sym_n}$ and $S(V)^{\text{WB}_n}$ needed in the following statements.)

%%%%%%%%%%%%%%%%%%%%%%%%%%%%%%%%%%%%%%%%%%%%%
\begin{cor}
Let $G$ be the symmetric group $\Sym_n=G(1,1,n)$ for $n\geq 4$
acting on $V$  with $\CC$-basis
$v_1, \ldots, v_n$ (via the natural permutation representation).
The Hochschild cohomology in degree $2$ for the skew group algebra
$S(V)\#G$ is
$$
  \HH^2(S(V)\#G) \ \cong
  \ \HH^2(1) \oplus \HH^2(1,2,3)
$$
where 
\begin{itemize}
\item
$
\HH^2(1) = 
     \left(S(V) \otimes  \Wedge^2 (V^*)\right)^{\Sym_n}
   = 
     \displaystyle\bigoplus_{1\leq i_1 < i_2 \leq n}\ 
     S(V)^{\Sym_n}\ (\theta_{i_1} \wedge\theta_{i_2}),
$\\
for
$$
 \theta_j := \sum_{1\leq i\leq n} v_i^{j-1} \otimes  x_i,
$$
\item
$
  \HH^2(1,2,3)=
  \ \CC[f_0,f_1,\ldots, f_{n-3}],
$
where
$f_0 = (v_1 + v_2 + v_3)$
and $f_i$ 
is the $i$-th elementary symmetric function in
$v_{4}, \ldots, v_{n}$ for $i>0$.
\end{itemize}
\end{cor}
%%%%%%%%%%%%%%%%%%%%%%%%%%%%%%%%%%%%%%%%%%%%%%%%%%%%%%%%%%%%%%%%%%%

%%%%%%%%%%
% Weyl Group B_n case:
%%%%%%%%%
%%%%%%%%%%%%%%%%%%%%%%%%%%%%%%%%%%%%%%%%%%%%%
\begin{cor}
Let $G$ be the Weyl group ${WB}_n=G(2,1,n)$ for $n\geq 4$
acting on $V$  with $\CC$-basis
$v_1, \ldots, v_n$ (via its natural reflection representation).
The Hochschild cohomology in degree $2$ for the skew group algebra
$S(V)\#G$ is 
$$
  \HH^2(S(V)\#G) \ \cong 
  \ \HH^2(1) \oplus \HH^2(1,2,3)\oplus \HH^2(1, -2) 
$$
where 
\begin{itemize}
\item
$
\HH^2(1) = 
     \left(S(V) \otimes  \Wedge^2 (V^*)\right)^{WB_n}
   = 
     \displaystyle\bigoplus_{1\leq i_1 < i_2 \leq n}\ 
     S(V)^{WB_n}\ (\theta_{i_1} \wedge\theta_{i_2}),
$\\
for
$$
 \theta_j := \sum_{1\leq i\leq n} (v_i)^{2(j-1)+1} \otimes  x_i,
$$
\item
$
  \HH^2(1,2,3)=
  \ \CC[f_0,f_1,\ldots, f_{n-3}],
$
where
$f_0 = (v_1 + v_2 + v_3)^2$
and $f_i$ is the $i$-th elementary symmetric function of 
$v_{4}^2, \ldots, v_{n}^2$ for $i>0$,
\item
$
  \HH^2(1, -2)\ = \CC[f'_1, \ldots, f'_{n-2}] ,
$
where
$f'_i$ is the $i$-th elementary symmetric function in 
$v_3^2, \ldots, v_n^2$.
\end{itemize}
\end{cor}
%%%%%%%%%%%%%%%%%%%%%%%%%%%%%%%%%%%%%%%%%%%

%%%%%%%%%%%%%%%%%%%%%%%%%%%%%%%%%%%%%%%%%%%%%%%%%%%%%%%%%%%%%%%%%%%%%%%
%%%%%%%%%%%%%%%%%%%%%%%%%%%%%%%%%%%%%%%%%%%%%%%%%%%%%%%%%%%%%%%%%%%%%%%
%%%%%%%%%%%%%%%%%%%%%%%%%%%%%%%%%%%%%%%%%%%%%%%%%%%%%%%%%%%%%%%%%%%%%%%
%%%%%%%%%%%%%%%%%%%%%%%%%%%%%%%%%%%%%%%%%%%%%%%%%%%%%%%%%%%%%%%%%%%%%%%
%%%%%%%%%%%%%%%%%%%%%%%%%%%%%%%%%%%%%%%%%%%%%%%%%%%%%%%%%%%%%%%%%%%%%%%
\section{Invariant theory
for certain centralizer subgroups 
of  $G(r,p,n)$}
\label{semi-invariantsofG(r,p,n)actingfaithfully}

\vspace{2ex}

In Section~\ref{faithfulHochschildcohomologyindeg2},
we gave the Hochschild cohomology of $S(V)\# G$ 
in degree 2 where $V$ is the natural reflection representation 
of $G=G(r,p,n)$.
The Hochschild cohomology can be expressed as the direct sum of 
certain subspaces, each subspace invariant under an action of a 
centralizer subgroup $Z(g)$ for some $g \in G$.
In this section, we determine these invariant subspaces.
We analyze the contribution $\HH^2(g)$
from the group elements $g \in G(r,p,n)$ with $\codim(V^g)=2$
and $\det g=1$. 
Recall that $v_1, \ldots, v_n$ is the given $\CC$-basis for $V$. 
%%%%%%%%%%%%%%%%%%%%%%%%%%%%%%%%%%%%%%%%%
\begin{prop}
\label{HH^2(g_b)}
Assume $n\geq 4$ and $G=G(r,p,n)$.  
Let $g=(1, 2, 3)$.
Then
$$
  S(V^g)^{Z(g)} = 
 \bigoplus_{\substack{
            0\leq i<r, \ %\\ 
            %\rule[0ex]{0ex}{1.5ex}%strut
            0\leq j <m\\
            \rule[0ex]{0ex}{1.5ex}%strut
            i \equiv 3 j r/p \mod r}}
    \CC[f_0^r, f_1, \ldots, f_{n'-1}, f_{n'}^m]\ \ \ 
    f_0^i\ f_{n'}^j\   ,
$$
where
\begin{itemize}
\item
$n'=n-3$,
\item
$m=p/3$ if $3$ divides $p$; $m=p$ otherwise,
\item
$f_0=v_1+v_2+v_3$,
\item
$f_{n'}=(v_4 \cdots v_{n})^{r/p}$,
\item
$f_i$ is the $i$-th elementary symmetric function of 
$v_{4}^r, \ldots, v_{n}^r$ for $1\leq i \leq n'-1$.
\end{itemize}
\vspace{1ex}
\noindent
Let $\chi$ be the Hochschild character for $g$ as defined in (\ref{S:hochchar}).
The set of Hochschild semi-invariants for $g$ is
$$ 
 \HH^2(g)=
 S(V^g)^{\chi}
 \ = \
 \bigoplus_{\substack{
            0\leq i<r, \  %\\ 
            %\rule[0ex]{0ex}{1.5ex}%strut
            0\leq j <m\\
            \rule[0ex]{0ex}{1.5ex}%strut
            i \equiv 2+3 j r/p \mod r}}
    \CC[f_0^r, f_1, \ldots, f_{n'-1}, f_{n'}^m]\ \ \ 
    f_0^i f_{n'}^j\   .
$$
\end{prop}
%%%%%%%%%%%%%%%%%%%%%%%%%%%%%%%%%%%%%%%%%%%%%%%%%%%%%
%%%%%%%%%%%%%%%%%%%%%%%%%%%%%%%%%%%%%%%%%%%%%%%%

%%%%%%%%%%%%%%%%%%%%%%%%%%%%%%%%%%%%%%%%%%%%%%%%
\begin{proof}
Let $v_0=f_0=v_1+v_2+v_3$.  Note
that $V^g= \CC\text{-span}\{ v_0,v_4, \ldots, v_{n}\}$
and 
$(V^g)^\perp \subset  \CC\text{-span}\{ v_1,v_2,v_3\}$.
Let 
$V_A= \CC\text{-span}\{ v_1,v_2,v_3\}$
and $V_B=\CC\text{-span}\{ v_4,\ldots,v_n\}$.
One may use (\ref{centralizer}) to verify that the
centralizer of $g$ in $G$ is
$$
  Z(g)=\{h_A \oplus h_B:
        h_A\in\langle \xi I_3, (1, 2, 3)\rangle < G(r,1,3),\
        h_B\in G(r,1,n')\} \cap G(r,p,n).
$$
The centralizer $Z(g)$ includes the subgroup
$$
  Z'=\{I_3 \oplus h_B: h_B\in G(r,p,n')\}.
$$
Since $\chi(h)=1$ for all $h\in Z'$,
$ 
  S(V^g)^{\chi} \subset S(V^g)^{Z'}.
$
And since $Z'\subset Z(g)$,
$
  S(V^g)^{Z(g)} \subset S(V^g)^{Z'}.
$
We begin by finding the invariants of $S(V^g)$
under $Z'$.

The group $Z'$ acts on $V_B$ as the reflection group $G(r,p,n')$
with invariants
$$
  S(V_B)^{G(r,p,n')}=\CC[v_4, \ldots, v_{n}]^{G(r,p,n')}=\CC[f_1,\ldots, f_{n'}]
$$
by Proposition~\ref{invariantsofG(r,p,n)}.
The group $Z'$ fixes $v_0$, so acts 
on $V^g$ as the direct sum of the trivial group with
$G(r,p,n')$. Thus,
$
 S(V^g)^{Z'}
  = \CC[v_0, v_4, \ldots, v_{n}]^{Z'}
  = \CC[f_0, f_1, \ldots, f_{n'}]
$
by Remark~\ref{directsumofreflectiongroups}.
Hence,
$$
 S(V^g)^{Z(g)} \subset \CC[f_0, f_1, \ldots, f_{n'}]
 \quad\text{and}\quad
  S(V^g)^{\chi} \subset
  \CC[f_0, f_1, \ldots, f_{n'}].
$$

The polynomials $f_1, \ldots, f_{n'-1}$ are invariant under $Z(g)$.
But the polynomials $f_0$ and $f_{n'}$ pick up a scalar when $Z(g)$ acts.
We introduce some characters to record and compare the actions
of $Z(g)$ on $f_0$ and on $f_{n'}$.

Define a linear character on $Z(g)$ 
$$
  \tau: Z(g) \rightarrow \CC
$$
by $h(f_0)=\tau(h) \, f_0$ for all $h$ in $Z(g)$.
Fix $h=h_A\oplus h_B\in Z(g)$
and note that $h_A = (\xi^k I_3) (1, 2, 3)^j$, for some $j$ and $k$.
Since $\det(h^\perp)=\xi^{2k}$
(as $(\dim V^g)^\perp=2$),
$
 \tau^2(h)
         = \xi^{2k}
         = \det(h^\perp)
         = \chi(h). 
$
Hence,
$\chi=\tau^2 $.

Recall that $\psi(M)$ is the product of the nonzero entries of $M$,
for any matrix $M$.  
Since $M \in G(r,p,n)$ implies that $\psi(M)^{r/p}=1$, we see that $\psi(h_B)^{r/p} = \psi(h_A)^{-r/p}$
(as $\psi(h)= \psi(h_A)\psi(h_B)$).
Observe that $\psi(h_A)
        = \xi^{3k}=\tau(h)^3$. 
Then 
$$
 h(f_{n'}) 
  =[\psi(h_B)(v_4\cdots v_{n})]^{r/p}
=\psi(h_B)^{r/p}\ f_{n'}
=\psi(h_A)^{-r/p}\ f_{n'}
=\tau(h)^{-3r/p}\ f_{n'}.   
$$

In summary, the action of $Z(g)$ on
$f_0, \ldots, f_n'$ is given by
\begin{equation}
\begin{aligned}
  h(f_k) = 
    \begin{cases}
             \quad\quad\tau(h)\  f_k &\text{ for } k = 0\\ 
           \quad\quad\quad\quad\ f_k &\text{ for } 0<k<n' \\
       \tau(h)^{-3r/p}\ f_{k} &\text{ for } k=n'
    \end{cases}\\
\end{aligned}
\end{equation}
for all $h\in Z(g)$.
In particular, for any $i,j$,
\begin{equation}
\label{summary}
\begin{aligned}
  h(f_0^i f_{n'}^j) &= \tau(h)^{i-3jr/p} \ f_0^i f_{n'}^j,
  \quad\quad\text{and}\\
  \chi(h) &= \tau^2(h).
\end{aligned}
\end{equation}

We now give an explicit description of 
$S(V^g)^{Z(g)}$ and $S(V^g)^{\chi}$
as subsets of  $\CC[f_0,\ldots, f_{n'}]$.
Fix some $s$ in $\CC[f_0, \ldots, f_{n'}]$
and write $s$ as a polynomial in $f_0$ and $f_{n'}$ with coefficients
in $f_1, \ldots, f_{n'-1}$:
$$
  s = \sum_{i,j \geq 0}
      \ a_{ij}\ f_0^i\ f_{n'}^j
      \quad\quad \text{  where  }\quad
      a_{ij} \in \CC[f_1, \ldots, f_{n'-1}]
      \subset S(V^g)^{Z(g)}.
$$
Then $s$ is invariant (respectively $\chi$-invariant)
under $Z(g)$ exactly when each
$f_0^i f_{n'}^j$ with $a_{ij}\neq 0$
is invariant (respectively $\chi$-invariant)
under $Z(g)$.
The monomial $f_0^i f_{n'}^j$ is invariant under $Z(g)$ exactly when
$$
  \tau(h)^{i-3jr/p} = 1
  \quad\quad\text{ for all }
  h \in Z(g)
$$
by Equation~(\ref{summary}).
Similarly, the monomial $f_0^i f_{n'}^j$ is $\chi$-invariant under $Z(g)$ exactly when
$$
\begin{aligned}
  \tau(h)^2\ f_0^i f_{n'}^j
  = \chi(h)\ f_0^i f_{n'}^j
  &= h (f_0^i f_{n'}^j)
  =  \tau(h)^{i-3jr/p}\ f_0^i f_{n'}^j, 
   \quad\quad\quad\text{i.e.,}\\
  \tau(h)^{2+3jr/p-i}&=1
   \qquad\qquad\qquad\text{ for all }
  h \in Z(g)\ .
\end{aligned}
$$
But the set $\{ \tau(h): h\in Z(g)\}$ depends on $p$ and $n$.

Note that $\tau(h)$ is a power of $\xi$ for each $h\in Z(g)$
and that $h'=\xi_1\xi_2\xi_3\xi_4^{-3}\in Z(g)$
with $\tau(h')=\xi$.
Hence 
$$
  \{ \tau(h): h\in Z(g)\} = \langle \xi \rangle.
$$
Thus the monomial $f_0^i f_{n'}^j$ is invariant under $Z(g)$ exactly when
$\xi^{i-3jr/p} = 1$:
$$ 
 f_0^i f_{n'}^j \in S(V^g)^{Z(g)}
 \quad\quad\text{ if and only if }\quad\quad
 i\equiv 3j r/p \mod r.
$$
In particular, $f_0^r$ is invariant under $Z(g)$.
Also, $f_{n'}^m$ is invariant under $Z(g)$
where $m=p/3$ if $p$ is divisible by $3$ and $m=p$ otherwise.
So
$$
 \CC[f_0^r, f_1, \ldots, f_{n'-1}, f_{n'}^m]
 \subset S(V^g)^{Z(g)}
 \subset
  \CC[f_0, f_1, \ldots, f_{n'-1}, f_{n'}].
$$
Thus, $S(V^g)^{Z(g)}$
is the free $\CC[f_0^r, f_1, \ldots, f_{n'-1}, f_{n'}^m]$-module
$$ 
 S(V^g)^{Z(g)}
 \ = \
 \bigoplus_{\substack{
            0\leq i<r, \ % \\ 
            %\rule[0ex]{0ex}{1.5ex}%strut
            0\leq j <m\\
            \rule[0ex]{0ex}{1.5ex}%strut
            i \equiv 3 j r/p \mod r}}
    \CC[f_0^r, f_1, \ldots, f_{n'-1}, f_{n'}^m]\ \ \ 
    f_0^i f_{n'}^j\   .
$$

Similarly, $f_0^i f_{n'}^j$ is $\chi$-invariant under $Z(g)$ exactly when
$ \xi^{2+3jr/p-i}=1$:
$$ 
 f_0^i f_{n'}^j \in S(V^g)^{\chi}
 \quad\quad\text{ if and only if }\quad\quad
 i\equiv 2+3j r/p \mod r.
$$
Thus, 
$S(V^g)^{\chi}$ 
is the free $\CC[f_0^r, f_1, \ldots, f_{n'-1}, f_{n'}^m]$-module
$$ 
 S(V^g)^{\chi}
 \ = \
 \bigoplus_{\substack{
            0\leq i<r, \ % \\ 
            %\rule[0ex]{0ex}{1.5ex}%strut
            0\leq j <m\\
            \rule[0ex]{0ex}{1.5ex}%strut
            i \equiv 2+3 j r/p \mod r}}
    \CC[f_0^r, f_1, \ldots, f_{n'-1}, f_{n'}^m]\ \ \ 
    f_0^i f_{n'}^j\   .
$$

\end{proof}
%%%%%%%%%%%%%%%%%%%%%%%%

\begin{prop}
\label{HH^2(g_c)}
Assume $n\geq 4$, $G=G(r,p,n)$.
Suppose 
$g=(1, -2 )$. Then
$$
  S(V^g)^{Z(g)} = 
 \begin{cases}
  \ \CC[f_1, \ldots, f_{n'-1}, f_{n'}^{\, p/2}]\ \ \
 &\text{if $p$ is even}\\
  \ \CC[f_1, \ldots, f_{n'-1}, f_{n'}^{\, p}]\ \ \
 &\text{if $p$ is odd}
 \end{cases}
$$
where
\begin{itemize}
\item
$n'=n-2$,
\item
$f_{n'}=(v_3 \ldots v_{n})^{r/p}$,
\item 
$f_i$ is the $i$-th elementary symmetric function of
$v_{3}^r, \ldots, v_{n}^r$ for $i<n'$.
\end{itemize}

\vspace{1ex}

\noindent
Let $\chi$ be the Hochschild character for $g$ as defined in
(\ref{S:hochchar}).
The set of Hochschild semi-invariants for $g$ is
$$ 
 \HH^2(g)=
 S(V^g)^{\chi}
 \ = \
 \bigoplus_{\substack{
            0\leq i<r \\ 
            \rule[0ex]{0ex}{1.5ex}%strut
            2 \equiv 2(-r/p) i \mod r}}
    \CC[f_1, \ldots, f_{n'-1}, f_{n'}^r]\ \ \ 
            f_{n'}^i\   .
$$
In particular, $\HH^2(g)$ is zero unless $r=2p$.
\end{prop}
%%%%%%%%%%%%%%%%%%%%%%%%%%%%%%%%%%%%%%%%%%%%%%%%%%%%%

%%%%%%%%%%%%%%%%%%%%%%%%%%%%%%%%%%%%%%%%%%%%%%%%
\begin{proof}
Each element $h$ of the centralizer $Z(g)$ breaks into a direct sum:
We may write $h=h^\perp \oplus h_\perp$,
where $h^\perp \in G(r,1,2)$ is the matrix of
$h$ acting on $(V^g)^{\perp}$ (with respect to the basis
$\{v_1, v_2 \}$ of $(V^g)^{\perp}$) 
and $h_\perp \in G(r,1,n')$ is the matrix of
$h$ acting on $V^g$ (with respect to the basis
$\{v_3, \ldots, v_n \}$ of $V^g$).
One may use (\ref{centralizer}) to verify that 
$$
  Z(g) = \{
          h^\perp\oplus h_\perp:
          h^\perp \in 
          \langle g^\perp, \xi I_2 \rangle,
          \ h_\perp\in G(r,1,n')
         \} \ \
         \cap\ \ G(r,p,n).
$$

The centralizer $Z(g)$ includes the subgroup $Z'$ of $G$ fixing $(V^g)^\perp$ 
pointwise:
$$
  Z'=\{I_2 \oplus B: B\in G(r,p,n')\} \ \subset\ Z(g).
$$
The group $Z'$ acts on $V^g=\CC\text{-span}\{v_3, \ldots, v_n\}$ 
as the reflection group $G(r,p,n')$.
By Proposition~\ref{invariantsofG(r,p,n)},
$$
  S(V^g)^{Z'} = \CC[v_3, \ldots, v_n]^{G(r,p,n')}
             =\CC[f_1, \ldots, f_{n'}].
$$
Since $\chi(h)=1$ for all $h\in Z'$,
$ 
  S(V^g)^{\chi} \subset S(V^g)^{Z'}.
$
And since $Z'\subset Z(g)$,
$
  S(V^g)^{Z(g)} \subset S(V^g)^{Z'}.
$
Hence,
$$
 S(V^g)^{Z(g)} \subset \CC[f_1, \ldots, f_{n'}]
 \quad\text{and}\quad 
  S(V^g)^{\chi} \subset
  \CC[f_1, \ldots, f_{n'}].
$$

Observe that
$f_1, \ldots, f_{n'-1}$ lie in $S(V^g)^{Z(g)}$.
How does $Z(g)$ act on the last polynomial, $f_{n'}$?
Recall that $\psi(M)$ is
the product of the nonzero entries of $M$, for any matrix $M$,
and that $M \in G(r,p,n)$ implies $\psi(M)^{r/p}=1$.
For any $h \in Z(g)$, $\psi(h)= \psi(h^\perp)\psi(h_\perp)$ and
\begin{equation}
\label{lastpoly-g_bcase}
\begin{aligned}
  h(f_{n'}) 
  &=h(v_3 \cdots v_n)^{r/p}\\
  &= (\psi(h_\perp) v_3 \cdots v_n)^{r/p}\\
  &=\psi(h_\perp)^{r/p}\ f_{n'}\\
  &= \psi(h)^{r/p}\, \psi(h^\perp)^{-r/p}\ f_{n'}\\
  &= \psi(h^\perp)^{-r/p}\ f_{n'}.
\end{aligned}
\end{equation}

Fix some $s$ in $\CC[f_1, \ldots, f_{n'}]$
and write $s$ as a polynomial in $f_{n'}$ with coefficients
in $f_1, \ldots, f_{n'-1}$:
$$
  s = \sum_{i \geq 0}
      \ a_{i}\ f_{n'}^i
      \quad\quad \text{  where  }\quad
      a_{i} \in \CC[f_1, \ldots, f_{n'-1}]
      \subset S(V^g)^{Z(g)}.
$$
Then $s$ is invariant (respectively $\chi$-invariant)
under $Z(g)$ exactly when each
$f_{n'}^i$ with $a_{i}\neq 0$
is invariant (respectively $\chi$-invariant)
under $Z(g)$.

From Equation~(\ref{lastpoly-g_bcase}),
the polynomial $f_{n'}^i$ is invariant under  $Z(g)$ exactly when
\begin{equation}
\label{1-lastpoly-g_bcase}
  1 = \psi(h^\perp)^{-ir/p}
  \quad\quad\text{ for all }
  h \in Z(g).
\end{equation}
And the polynomial $f_{n'}^i$ is $\chi$-invariant under $Z(g)$ exactly when
\begin{equation}
\begin{aligned}
  \det(h^\perp)\ f_{n'}^i
  &= \chi(h)\ f_{n'}^i
  = h (f_{n'}^i)
  =  \psi(h^\perp)^{-ir/p}\ f_{n'}^i, \quad\quad\text{ i.e.,}\\
\label{2-lastpoly-g_bcase}
  \det(h^\perp) &= \psi(h^\perp)^{-ir/p},
   \qquad\text{ for all }
  h \in Z(g)\ .
\end{aligned}
\end{equation}

Note that $Z(g)^\perp:=\{h^\perp: h \in Z(g)\}$
is generated by $g^\perp$ and $\xi I_2=\xi_1\xi_2$
since $n\geq 3$ implies that the centralizer $Z(g)$ contains
$\xi_1\xi_2 \xi_3^{-2} \in G(r,p,n)$.
(If $n=2$, the matrix $\xi I_2$ need not lie in $G(r,p,n)$.)
Thus, in order to establish the invariance or semi-invariance of $f_{n'}^i$,
we need only check Equations~(\ref{1-lastpoly-g_bcase}) 
and (\ref{2-lastpoly-g_bcase}) for $g^\perp$ and $\xi I_2$
Note that $g = \xi_2^{r/2}(1, 2)\in G(r,p,n)$ implies that the integer
$r/p$ is even.
Hence $\psi(g^{\perp})^{-ir/p}=(-1)^{-ir/p}=1$.
And $\det(g^\perp)=1$.
Hence, both  Equations~(\ref{1-lastpoly-g_bcase}) 
and (\ref{2-lastpoly-g_bcase}) are satisfied for the 
generator $g^\perp$ of $Z(g)^\perp$.
So, the invariance or semi-invariance of 
$f_{n'}^i$ depends only
on $i$ and the action of the generator $\xi I_2$ of $Z(g)^{\perp}$
on $f_{n'}$.

Thus, $f_{n'}^i$ is invariant under $Z(g)$ exactly when
$$
  1  =\psi(\xi I_2)^{-ir/p} 
  =(\xi^2)^{-ir/p}, 
   \qquad\text{i.e.}\quad
  0 \equiv 2 i r/p \mod r.
$$
But $0 \equiv 2 i r/p \mod r$ means that $i$ is a multiple of $p/2$ 
when $p$ is even and a multiple of $p$ when $p$ is odd.
And $f_{n'}^i$ is $\chi$-invariant under $Z(g)$ exactly when
$$
  \xi^2
  =\det(\xi I_2)
  =\psi(\xi I_2)^{-ir/p} 
  =(\xi^2)^{-ir/p}, 
   \qquad\text{i.e.}\quad
  2\equiv -2 i r/p \mod r.
$$
In particular, $f_{n'}^r$ is invariant under $Z(g)$.
Hence
$$
  S(V^g)^{Z(g)} = 
 \begin{cases}
  \ \CC[f_1, \ldots, f_{n'-1}, f_{n'}^{\, p/2}]\ \ \
 &\text{if $p$ is even}\\
  \ \CC[f_1, \ldots, f_{n'-1}, f_{n'}^{\, p}]\ \ \
 &\text{if $p$ is odd}
 \end{cases}
$$
and
$$ 
 \HH^2(g)=
 S(V^g)^{\chi}
 \ = \
 \bigoplus_{\substack{
            0\leq i<r \\ 
            \rule[0ex]{0ex}{1.5ex}%strut
            2 \equiv 2(-r/p) i \mod r}}
    \CC[f_1, \ldots, f_{n'-1}, f_{n'}^r]\ \ \ 
            f_{n'}^i\   .
$$
\end{proof}
%%%%%%%%%%%%%%%%%%%%%%%%%

%%%%%%%%%%%%%%%%%%%%%%%%%%%%%%%%%%%%%%%%%
\begin{prop}
\label{HH^2(g_d)1}
Assume $n\geq 2$, $r\geq 2$, $G=G(r,p,n)$.
Let $g=\xi_1^{\ell_1} \xi_2^{\ell_2}$, where $\ell_i \neq 0$.
The set of Hochschild semi-invariants for $g$ is trivial,
$$ 
\HH^2(g)= S(V^g)^\chi=0,
$$
unless $G=G(r,r,n)$.
\end{prop}
%%%%%%%%%%%%%%%%%%%%%%%%%%%%%%%%
\begin{proof}
Note that $V^g = \CC\text{-span}\{v_3, \ldots, v_n\}$ and
 $(V^g)^\perp=\CC\text{-span}\{v_1, v_2\}$.
Let $h = \xi_1 \xi_2^{p-1} \in Z(g)$.
Then $h|_{V^g}=1$
and $\det h^{\perp}=\xi^p$.
Lemma~\ref{forceszero} implies that if $\HH^2(g)$ is nonzero, then
$\xi^p=1$ and so $p=r$.
\end{proof}

%%%%%%%%%%%%%%%%%%%%%%%%%%%%%%%%%%%%%%%%%
\begin{prop}
\label{HH^2(g_d)2}
Assume $n\geq 2$, $r\geq 2$ ($r$ even).
Let $G=G(r,r,n)$.
Let $g=\xi_1^{r/2} \xi_2^{r/2}$.
The set of Hochschild semi-invariants for $g$ is trivial:
$$ 
\HH^2(g)= S(V^g)^\chi=0.
$$
\end{prop}
%%%%%%%%%%%%%%%%%%%%
\begin{proof}
Note that $V^g = \CC\text{-span}\{v_3, \ldots, v_n\}$ and
$(V^g)^\perp=\CC\text{-span}\{v_1, v_2\}$.
Let $h= (1, 2)$.
Then $h \in Z(g)$ with $h|_{V^g}=1$.
But $\det h^\perp =-1 \neq 1$.
By Lemma~\ref{forceszero},
$\HH^2(g)$ is zero.
\end{proof}

%%%%%%%%%%%%%%%%%%%%%%%%%%%%%%%%%%%%%%%%%
\begin{prop}
\label{HH^2(g_d)3}
Assume $n\geq 3$, $r\geq 2$.
Let $G=G(r,r,n)$.
Let $g=\xi_1^{\ell} \xi_2^{-\ell}$, 
where $\ell \neq r/2$.
Then
$$
 S(V^g)^{Z(g)} = \CC[f_1, \ldots, f_{n'-1}, f_{n'}^r],
$$
where
\begin{itemize}
\item
$n'=n-2$,
\item
$f_{n'}=v_3 \cdots v_{n}$,
\item
$f_i$ is the $i$-th elementary symmetric function
of $v_3^r, \ldots, v_{n}^r$ for $i<n'$.
\end{itemize}
Let $\chi$ be the Hochschild character for $g$ as defined in (\ref{S:hochchar}).
The set of Hochschild semi-invariants for $g$ is
$$ 
 \HH^2(g)=
 S(V^g)^{\chi}
 \ = \
    \CC[f_1, \ldots, f_{n'-1}, f_{n'}^r]\ \ \ 
    f_{n'}^{r-1}\   .
$$
\end{prop}
%%%%%%%%%%%%%%%%%%%%
\begin{proof}
Note that $V^g = \CC\text{-span}\{v_3, \ldots, v_n\}$ and
$(V^g)^\perp=\CC\text{-span}\{v_1, v_2\}$.
Let $s\in \HH^2(g)$ be a Hochschild semi-invariant for $g$.
Consider the subgroup of $Z(g)$
$$
  Z'=\{I_2 \oplus B: B\in G(r,r,n')\}.
$$
Since $\det h^\perp=1$ for all $h\in Z'$,
$s$ is invariant under $Z'$: 
$s\in S(V^g)^{Z'}=\CC[v_3,\ldots, v_n]^{Z'}$.
But $Z'$ acts on the vector space $V^g$ as $G(r,r,n')$
with respect to the basis $\{v_3, \ldots, v_n\}$ of $V^g$.
By Proposition~\ref{invariantsofG(r,p,n)}, 
$$
  S(V^g)^{Z'}= \CC[v_3,\ldots, v_n]^{G(r,r,n')} = \CC[f_1, \ldots, f_{n'}].
$$
Hence,
$$
 \HH^2(g)=\{ \text{Hochschild semi-invariants for } g \}
 \subset
  \CC[f_1, \ldots, f_{n'}].
$$

One may use (\ref{centralizer}) to verify that
$$
 Z(g)=\{A\oplus B\in G(r,r,n): 
       A\text{ is a } 2\times 2\text{ diagonal matrix} \}.
$$
Let $s$ be any element of $\CC[f_1, \ldots, f_{n'}]$ and write
$s$ as a polynomial in $f_{n'}$ with coefficients
in $f_1, \ldots, f_{n'-1}$:
$$ 
 s = \sum_{i\geq 0}\ a_i\ f_{n'}^i 
 \qquad\qquad\text{ where each }
  a_i \in \CC[f_1, \ldots, f_{n'-1}].
$$
Note that $f_1,\ldots, f_{n'-1}$ are invariant under $Z(g)$.
Hence $s$ is invariant (respectively $\chi$-invariant)
under $Z(g)$ exactly when each
$f_{n'}^i$ with $a_{i}\neq 0$
is invariant (respectively $\chi$-invariant)
under $Z(g)$.

Observe that 
$$
h(f_{n'}) = \chi(h)^{-1} f_{n'}= (\det h^\perp)^{-1}\ f_{n'}
$$
for any $h\in Z(g)$ since $p=r$.
Hence,
$f^i_{n'}$ is invariant under $Z(g)$ exactly when 
$$
 1 = (\det h^{\perp})^{-i} \quad\text{ for all } h \in Z(g).
$$
Similarly,
$f_{n'}$ is $\chi$-invariant under $Z(g)$ exactly when 
$$
\begin{aligned}
 \det(h^\perp) &= (\det h^{\perp})^{-i} \quad\text{ for all } h \in Z(g),
 \qquad\qquad\text{i.e.,}\\
 1 &= (\det h^{\perp})^{1+i} \quad\text{ for all } h \in Z(g).
\end{aligned}
$$
As $n\geq 3$,
$\{ \det h^\perp : h \in Z(g) \} = \langle \xi \rangle.
$
Hence, $f_{n'}^i$ is invariant under $Z(g)$ exactly when 
$i$ is a multiple of $r$.
And $f_{n'}^i$ is $\chi$-invariant under $Z(g)$ exactly when 
$  -1 \equiv i \mod r $.
Thus,
$$
 S(V^g)^{Z(g)} = \CC[f_1, \ldots, f_{n'-1}, f_{n'}^r],
$$
and
$$ 
 \HH^2(g)=
 S(V^g)^{\chi}
 \ = \
    \CC[f_1, \ldots, f_{n'-1}, f_{n'}^r]\ \ \ 
    f_{n'}^{r-1}\   .
$$
\end{proof}
%%%%%%%%%%%%%%%%%%%%%%%%%

%%%%%%%%%%%%%%%%%%%%
\begin{prop}
\label{HH^2(g_e)}
Assume $r\geq 2$, $G=G(r,p,n)$.
Let $g=(1, 2) \xi_3^{\ell}$.
The set of Hochschild semi-invariants for $g$ is zero:
$$ 
\HH^2(g)= S(V^g)^\chi=0.
$$
\end{prop}
%%%%%%%%%%%%%%%%%%%%%%%%%%%%%%%%%%%%%%%%%%%%%%%%
\begin{proof}
Note that $V^g = \CC\text{-span}\{v_1+v_2, v_4, v_5, \ldots, v_n\}$ and
$(V^g)^\perp=\CC\text{-span}\{v_1-v_2, v_3\}$.
Let $h=(1, 2) \in Z(g)$.
Then $h|_{V^g}=1$ while $\det h^\perp = -1$,
so Lemma~\ref{forceszero}
implies that $\HH^2(g)=0$.
\end{proof}
% end type e
%%%%%%%%%%%%%%%%%%%%%%%%%

%%%%%%%%%%%%%%%%%%%%
\begin{prop}
\label{HH^2(g_f)}
Assume $n\geq 4$.
Let $g=(1, 2) \xi_3^{a}\xi_4^{-a}(3, 4)$.
The set of Hochschild semi-invariants for $g$ is zero:
$$ 
\HH^2(g)= S(V^g)^\chi=0.
$$
\end{prop}
%%%%%%%%%%%%%%%%%%%%%%%%%%%%%%%%%%%%%%%%%%%%%%%%
\begin{proof}
Note that $V^g = \CC\text{-span}\{v_1+v_2, v_3 + \xi^{-a} v_4, 
v_5, \ldots, v_n\}$ and
that $(V^g)^\perp=\CC\text{-span}\{v_2-v_1, v_3-\xi^{-a} v_4\}$.
Let $h=(1, 2) \in Z(g)$.
Then $h|_{V^g}=1$ while $\det h^\perp = -1$,
so Lemma~\ref{forceszero}
implies that $\HH^2(g)=0$.
\end{proof}
% end type f
%%%%%%%%%%%%%%%%%%%%%%%%%  

%%%%%%%%%%%%%%%%%%%%%%%%%%%%%%%%%%%%%%%%%%%%%%%%%%%%%%%%%%%%%%%%%%%%%%%
%%%%%%%%%%%%%%%%%%%%%%%%%%%%%%%%%%%%%%%%%%%%%%%%%%%%%%%%%%%%%%%%%%%%%%%
%%%%%%%%%%%%%%%%%%%%%%%%%%%%%%%%%%%%%%%%%%%%%%%%%%%%%%%%%%%%%%%%%%%%%%%
%%%%%%%%%%%%%%%%%%%%%%%%%%%%%%%%%%%%%%%%%%%%%%%%%%%%%%%%%%%%%%%%%%%%%%%
%%%%%%%%%%%%%%%%%%%%%%%%%%%%%%%%%%%%%%%%%%%%%%%%%%%%%%%%%%%%%%%%%%%%%%%

\section{Hochschild 2-cohomology for nonfaithful action 
of $G(r,p,n)$}\label{hh2nonfaithful}
Let $G=G(r,1,n)$, with $n\geq 3$.  We define a {\em nonfaithful} action of $G$
on $V=\CC^n$ and 
determine the Hochschild $2$-cohomology for the skew group algebra
$S(V) \# G$. This nonfaithful action will yield nontrivial graded
Hecke algebras in Section \ref{ghaG(r,1,n)}.

Although we restrict our attention to $G(r,1,n)$ in this section, the same results generally hold for $G(r,p,n)$ (under the nonfaithful action).
In fact, Propositions \ref{diagonalmatrixcontributation-nonfaithfulcohomology} and
\ref{nondiagonalmatrixcontributation-nonfaithfulcohomology}
and Theorem \ref{main-nonfaithful-theorem}
  below hold verbatim when we replace $G(r,1,n)$ by $G(r,p,n)$ and replace
$n\geq 3$ by $n\geq 5$.
We simply intersect each centralizer $Z(g)$ in $G(r,1,n)$ with $G(r,p,n)$ throughout the proofs.
We leave the details to the reader.
(For $p\geq 2$ and low values of $n$, Remark \ref{diagonalmatrixcontributation-nonfaithfulcohomology-remark}
  does not apply as the centralizers exhibit a different structure.)

Recall that $G$ is the group product $N\cdot G(1,1,n)$,
where $N$ is the normal subgroup of diagonal matrices in $G$
and $G(1,1,n)\leq G$ is isomorphic to $\Sym_n$.
Define a representation 
\begin{equation}\label{S:rep}
\rho: G\rightarrow {\rm GL}(V)
\end{equation}
by composing
the quotient map $G\rightarrow \Sym_n$ with the permutation
representation of $\Sym_n$ on $V$ (permuting the fixed basis $v_1,\ldots,v_n$).

The next lemma describes the centralizer of a diagonal element.
Given a diagonal matrix
$g=\xi^{a_1}I_{n_1}\oplus\cdots\oplus\, \xi^{a_k}I_{n_k}$,
rename
the basis vectors $v_1, \ldots, v_n$ of $V$ 
so that $v_1^{(1)}$, \ldots, $v_{n_1}^{(1)}$ are the first $n_1$ basis vectors,
$v_1^{(2)}, \ldots,v_{n_2}^{(2)}$ are the next $n_2$ basis vectors, etc.,
and
decompose $V$ as $V_1\oplus\cdots\oplus V_k$ with 
$V_i = \CC\text{-span}\{v_1^{(i)},\cdots,v_{n_i}^{(i)}\}$.
One may use (\ref{centralizer}) to verify the following lemma.
%%%%%%%%%%%%%%%%%%%%%%%%%%%%%%%%%%%%%%%%%%%%%%%%%
\begin{lemma}
\label{determinationofZ(g)}
Let $g'$ be a diagonal matrix in $G=G(r,1,n)$
and suppose $n_1, \ldots, n_k$ are the multiplicities
of the diagonal entries.
Then $g'$ is conjugate in $G$ to some 
$g= \xi^{a_1} I_{n_1} \oplus\cdots\oplus \xi^{a_k} I_{n_k}$ and
$Z(g)$ acts on $V=V_1\oplus \cdots\oplus V_k$ as the direct sum
of reflection groups,
$$
Z(g)  = \{ M_1 \oplus\cdots\oplus M_k: 
           \  M_i \in G(r,1,n_i)
        \} %\\
      = G(r,1,n_1) \oplus \cdots \oplus G(r,1,n_k).
$$
\end{lemma}
%%%%%%%%%%%%%%%%%%%%%%%%%%%%%%%%%%%%%%%%%%%%%%%%%%%

%%%%%%%%%%%%%%%%%%%%%%%%%%%%%%%%%%%%%%%%%%%%%%%%%%%%%%%%%%%%%%%%%%%%%
%%%%%%%%%%%%%%%%%%%%%%%%%%%%%%%%%%%%%%%%%%%%%%%%%%%%%%%%%%%%%%%%%%%%%
\subsection*{Contribution from $g$ with $\codim V^g=0$}
\label{diagonalmatrixcontributation-nonfaithfulcohomology-subsection}

\rule[0ex]{0ex}{1ex}%strut

We first find the contribution to the Hochschild $2$-cohomology
of $S(V)\# G(r,1,n)$
from $g$ with $\codim V^g=0$.
The set of $g$ in $G=G(r,1,n)$ with $\codim(V^g)=0$ under
the nonfaithful representation $\rho$ is exactly the
subgroup $N$ of diagonal matrices in $G(r,1,n)$.
Below, we use the notation $\HH^2(g)$ defined
in (\ref{S:hh2g}) and we identity
$S(V)\otimes\wedge^k (V^*)$ with the $S(V)$-module of $k$-derivations
on $V^*$.
%%%%%%%%%%%%%%%%%%%%%%%%%%%%%%%%%%
\begin{prop}
 \label{diagonalmatrixcontributation-nonfaithfulcohomology}
Let $G=G(r,1,n)$ act on $V = \CC^n$ via the nonfaithful
representation $\rho$.  Let $g'\in G$ be a diagonal matrix.
The $g'$ is conjugate to some $g\in G$ whose set of Hochschild semi-invariants
$\HH^2(g)$ is the set of derivation $2$-forms invariant under a direct sum of symmetric groups.

Explicitly: Suppose $n_1, \ldots, n_k$ are the multiplicities of the diagonal entries of $g'$ and
decompose $V$ as $V_1\oplus\cdots\oplus V_k$ with each $V_i \cong \CC^{n_i}$ 
as above. 
For $i=1,\ldots,k$, 
let $ \theta^{(i)}_{1}, \ldots, \theta^{(i)}_{n_i} \in S(V_i)\otimes V_i^*$
be basic derivations for the permutation action of the symmetric group
$\Sym_{n_i}$ on $V_i$.
Then $\HH^2(g)$ is generated by the wedge products of the $\theta_j^{(i)}$ taken two at a time:
$$
\begin{aligned}
  \HH^2(g)\ = \ 
  \left( S(V) \otimes  \Wedge^2 (V^*) \right)^{Z(g)}
   &=
  \bigoplus_{\substack{
i \leq i', \\ j < j' \text{ when } i=i'}}
  S(V)^{Z(g)}\  (\theta^{(i)}_{j}\wedge  \theta^{(i')}_{j'})\, . 
\end{aligned}
$$
Also,
$$
  S(V)^{Z(g)} = \bigoplus_{i=1}^{k}\, S(V_i)^{\Sym_{n_i}}.
$$
\end{prop}
%%%%%%%%%%%%%%%%%%%%%%%%%%%%%%%%%%
\begin{remark}
\label{diagonalmatrixcontributation-nonfaithfulcohomology-remark}
We may easily construct explicit basic invariants and basic derivations for the last proposition using Propositions \ref{invariantsofG(r,p,n)} and \ref{explicitderivations}.
Let $x_1^{(i)}, \ldots, x_{n_i}^{(i)}$ be the basis of $V_i^*$
dual to the basis $v_1^{(i)},\ldots, v_{n_i}^{(i)}$ 
of $V_i$.
Let $f^{(i)}_j$  be the $j$-th elementary symmetric
function of $v_1^{(i)},\ldots ,v_{n_i}^{(i)}$.
Then
$$
  S(V)^{Z(g)} = \CC[f^{(1)}_1, \ldots, f^{(1)}_{n_1}, \ldots,
                 f^{(k)}_1, \ldots, f^{(k)}_{n_k}].
$$
We may also explicitly define 
$\theta_1^{(i)}, \ldots, \theta_{n_i}^{(i)} \in S(V_i)\otimes V_i^*$:
$$
 \theta^{(i)}_{j} := 
  \sum_{1\leq l\leq n_i} (v_l^{(i)})^{(j-1)r+1} \otimes  x_l^{(i)}.
$$
\end{remark}
%%%%%%%%%%%%%%%%%%%%%%%%%%%%%%%%%%%%%%%%%%%%%%%%%%%
%
\begin{proof}[Proof of Proposition
\ref{diagonalmatrixcontributation-nonfaithfulcohomology}]
By Lemma~\ref{determinationofZ(g)}, 
$g'$ is conjugate to some diagonal $g$ 
with
$$
 Z(g)=G(r,1,n_1)\oplus\cdots\oplus G(r,1,n_k).
$$ 
Under the representation $\rho$, $g$ acts as the identity and $V^g=V$.  Thus, the Hochschild 
character $\chi: Z(g) \rightarrow \CC$ (given by $h\mapsto \det h^\perp$) is trivial and
$\HH^2(g)$ is the just the set of derivation $2$-forms invariant under $Z(g)$,
$$
  \HH^2(g)= 
  \left( S(V^g) \otimes \Wedge^{2 - \codim V^g} ((V^g)^*) \right)^{\chi}=
\left( S(V) \otimes  \Wedge^2 (V^*) \right)^{Z(g)} \ .
$$

Each $G(r,1,n_i)$ acts on $V_i$ as the symmetric group $\Sym_{n_i}\cong G(1,1,n_i)$ under the representation $\rho$.
Thus the block diagonal group $Z(g)$ acts as the direct sum of symmetric groups,
and we simply find the polynomials and derivation $2$-forms invariant
under the permutation action of
$
\Sym_{n_1} \oplus\cdots\oplus\Sym_{n_k}
$
on $V_1\oplus\cdots\oplus V_{k}$.

By Remark~\ref{directsumofreflectiongroups}, the ring of invariant polynomials is just the tensor product of the corresponding rings $S(V_i)^{\Sym_{n_i}}$ of symmetric polynomials.
The underlying arrangement of reflecting hyperplanes 
is given by the product of the polynomials defining the subarrangements.  More precisely, we define
the reflection arrangement 
(after identifying $(V^*)^*$ with $V$)
by the polynomial 
$$
  Q = Q_1 \cdots Q_k \ \in S(V),
$$
where each $Q_i\in S(V_i)$ 
is the Vandermonde determinant
($Q_i = \prod_{j_1<j_2} x^{(i)}_{j_1} - x^{(i)}_{j_2} $) 
for $\Sym_{n_i}$ acting on $V^*_i$ 
(see Section \ref{invthyofG(r,p,n)}).
As the $\theta_j^{(i)}$ are basic derivations,
the coefficient matrix of 
$\{ \theta_1^{(i)}, \ldots, \theta_{n_i}^{(i)}\}$
has determinant $Q_i$ up to a nonzero scalar
(by Theorem~\ref{Solomon'sThm}).
Thus
$
 \theta_1^{(1)}\wedge\cdots\wedge \theta_{n_k}^{(k)}$ is a nonzero scalar
 multiple of
 $(Q_1\cdots Q_k)\  x_1 \wedge\cdots\wedge x_{n}
 = Q\ x_1 \wedge\cdots\wedge x_{n}$.
By Proposition~\ref{Solomon'sThm}, the derivations $\theta^{(i)}_j$ wedged together two at a time then generate the invariant derivation $2$-forms over the ring of invariant polynomials. 
\end{proof}
%%%%%%%%%%%%%%%%%%%%%%%%%

%%%%%%%%%%%%%%%%%%%%%%%%%%%%%%%%%%%%%%%%%%%%%%%%%%%%%%%%%%%%%%%%%%%%%%%%%
%%%%%%%%%%%%%%%%%%%%%%%%%%%%%%%%%%%%%%%%%%%%%%%%%%%%%%%%%%%%%%%%%%%%%%%%%
\subsection*{Contribution from $g$ with $\codim V^g=2$}
\label{nondiagonalmatrixcontributation-nonfaithfulcohomology-subsection}

\rule[0ex]{0ex}{1ex}%strut

We now determine the contribution 
to the Hochschild $2$-cohomology of the skew group algebra
$S(V)\# G(r,1,n)$ (under the nonfaithful representation $\rho$
defined in (\ref{S:rep})) from $g$ with $\codim V^g=2$. 

%%%%%%%%%%%%%%%%%%%%%%%%%%%%%%%%%%%%%%%%%%%%%%%%%%%%%%%%%%%%%%
\begin{remark}
Suppose $g \in G=G(r,1,n)$ with $\codim(V^g)=2$
under the nonfaithful action of $G$ on $V$ given by $\rho$.
Then for $n\geq 3$, $g$ is conjugate to a diagonal matrix times
either a $3$-cycle or the product of two $2$-cycles (in case $n\geq 4$).  
In fact, we may assume $g$ is conjugate to a diagonal matrix times
the permutation
$$
  (1, 2, 3) \qquad\text{ or }\qquad (1, 2)(3, 4).
$$
If $g$ is the product of a diagonal matrix with  $(1, 2)(3, 4)$,
then the set of Hochschild semi-invariants for $g$ is zero,
i.e.\ $\HH^2(g)=0$:
Suppose $g=\xi_1^{a_1}\cdots \xi_n^{a_n}(1,2)(3,4)$ and define
$h=\xi_1^{a_1}\xi_2^{a_2}(1,2)\in Z(g)$.
Then $\det (h^{\perp})=-1 \neq 1$
(under the nonfaithful action given by $\rho$) and hence
$\HH^2(g)=0$ by Lemma~\ref{forceszero}.
\end{remark}
%%%%%%%%%%%%%%%%%%%%%%%%%%%%%%%%%%%%%%%%%%%%%%%%%%%%%%%%%%%%%%%

In the following proposition, we give explicit invariants by
decomposing $V$ into subspaces.
Let $V_A=\CC\text{-span}\{v_1,v_2,v_3\}$ and
$V_B=\CC\text{-span}\{v_4,\ldots,v_n\}$.
If $g'$ is the product of a diagonal matrix and a 3-cycle, then $g'$ is
conjugate to $g=g_A\oplus g_B$ where $g_A\in G(r,1,3)$ is 
a $3\times 3$ diagonal matrix times the 3-cycle $(1,2,3)$
and $g_B=\xi^{a_1}I_{n_1}\oplus\cdots \oplus \xi^{a_k}I_{n_k}
\in G(r,1,n-3)$.
(Note that $n_1,\ldots,n_k$ are simply the multiplicities of the diagonal entries of $g_B$.)
We apply Lemma~\ref{determinationofZ(g)} to $g_B$ acting on $V_B$:
Rename
the basis vectors $v_4, \ldots, v_n$ of $V_B$ 
so that $v_1^{(1)}$, \ldots, $v_{n_1}^{(1)}$ are the first $n_1$ basis vectors,
$v_1^{(2)}, \ldots,v_{n_2}^{(2)}$ are the next $n_2$ basis vectors, etc.;
decompose $V_B$ as $V_1\oplus\cdots\oplus V_k$ with 
$V_i = \CC\text{-span}\{v_1^{(i)},\cdots,v_{n_i}^{(i)}\}$.

%%%%%%%%%%%%%%%%%%%%%%%%%%%%%%%%%%%%%%%%%%%%%%%%%%%%%%%%%%%
\begin{prop}
\label{nondiagonalmatrixcontributation-nonfaithfulcohomology}
Assume $n\geq 3$.
Let $G=G(r,1,n)$ act on $V\cong \CC^n$ via the nonfaithful
representation $\rho$, and let $g'$ be the product of a diagonal matrix
with a $3$-cycle.
Let $f^{(i)}_l$ be the $l$-th elementary symmetric function
of the elements $v_j^{(i)}$ defined above.
Let $f_0=v_1+v_2+v_3$.
Then $g'$ is conjugate to $g$ with
$$
\begin{aligned}
  \HH^2(g)=S(V^g)^\chi
  = S(V^g)^{Z(g)} 
 &= \CC[f_0, f^{(1)}_1, \ldots, f^{(1)}_{n_1}, \ldots,
                 f^{(k)}_1, \ldots, f^{(k)}_{n_k}].
\end{aligned}
$$
\end{prop}
%%%%%%%%%%%%%%%%%%%%%%%%%%%%%%%%%%%%%%%%%%%%%%%%%%%%%%%%%%%%%%%%%
\begin{proof}
As in the above paragraph, let $g=g_A\oplus g_B$.
The action of $Z(g)$ on $V=V_A\oplus V_B$
(and on $V^g=V_A^{g} \oplus V_B^{g}$) decomposes as a direct sum:
$$
  Z(g)=Z_A \oplus Z_B
$$
where $Z_A=Z_{G(r,1,3)}(g_A)$ and $Z_B=Z_{G(r,1,n-3)}(g_B)$. 
In fact, 
$
  Z_A = \langle \xi I_3, (1,\, 2,\, 3) \rangle.
$
Under the representation $\rho$, $g$ acts on $V$ as the permuation
$(1,\, 2,\, 3)$.  Hence
$V^g=\CC\text{-span}\{v_1+v_2+v_3, v_4, \ldots, v_n\}$
and 
$(V^g)^\perp = \CC\text{-span}\{v_1-v_2, v_1-v_3\}$.
The subgroup $Z(g)$ acts on $(V^g)^\perp$
as the group $\langle (1,\, 2, \, 3) \rangle$ with determinant $1$ under the representation 
$\rho$.
Thus, the Hochschild 
character $\chi: Z(g) \rightarrow \CC$ (given by $h\mapsto \det h^\perp$) is trivial (the diagonal matrices act trivially):
$$
S(V^g)^\chi = S(V^g)^{Z(g)}.
$$

Note that $S(V_A^g)=\CC[v_1+v_2+v_3]$ and
$Z_{A}$ acts trivially on $v_1+ v_2+v_3$.
Hence, 
$$
   S(V_A^g)^{Z_{A}} =\CC[v_1+v_2+v_3].
$$ 
We apply
Lemma~\ref{determinationofZ(g)} 
and Remark~\ref{diagonalmatrixcontributation-nonfaithfulcohomology-remark}
to the diagonal matrix $g_B\in G(r,1,n-3)$ acting on $V_B$:
$$
  S(V_B^g)^{Z_B} = \CC[f^{(1)}_1, \ldots, f^{(1)}_{n_1}, \ldots,
                 f^{(k)}_1, \ldots, f^{(k)}_{n_k}].
$$
Hence, by Remark~\ref{directsumofreflectiongroups},
$$
  S(V^g)^{Z(g)} 
  \cong S( V_A^g)^{Z_A}
     \otimes
     S(V_B^g)^{Z_B}
 \cong \CC[f_0, f^{(1)}_1, \ldots, f^{(1)}_{n_1}, \ldots,
                 f^{(k)}_1, \ldots, f^{(k)}_{n_k}].
$$
\end{proof}
%%%%%%%%%%%%%%%%%%%%%%%%%%%%%%%%%%%%%%%%%%%%%%%%%%%%%%%%%

%%%%%%%%%%%%%%%%%%%%%%%%%%%%%%%%%%%%%%%%%%%%%%%%%%%%%%%%%%%%%%%%%%%%%%%%%
%%%%%%%%%%%%%%%%%%%%%%%%%%%%%%%%%%%%%%%%%%%%%%%%%%%%%%%%%%%%%%%%%%%%%%%%%
\subsection*{Putting it together, nonfaithful case}

\rule[0ex]{0ex}{1ex}%strut

We now record the Hochschild $2$-cohomology
of $S(V)\# G$ where $G=G(r,1,n)$ acts on $V$ 
via the nonfaithful representation $\rho$.
We combine the contribution to cohomology
from $g\in G$ with $\codim V^g=0$
and the contribution from $g\in G$ with $\codim V^g=2$.

%%%%%%%%%%%%%%%%%%%%%%%%%%%%%%%%%%%%%%%%%%%%%
\begin{theorem}\label{main-nonfaithful-theorem}
Assume $n\geq 3$.
Let $G=G(r,1,n)$ act on $V = \CC^n$ via the nonfaithful
representation $\rho$.
The Hochschild cohomology in degree $2$ for the skew group algebra
$S(V)\#G$ is 
$$
  \HH^2(S(V)\#G) \ \cong
  \ 
   \bigoplus_{g \in C}
       \HH^2(g) 
   \oplus 
   \bigoplus_{g' \in C'}
    \HH^2(g')
$$
where
\begin{itemize}
\item
  $C$ is a set of representatives of those
  conjugacy classes in $G$ containing diagonal matrices,
  and
  $\HH^2(g)$ is given in 
  Proposition~\ref{diagonalmatrixcontributation-nonfaithfulcohomology},
\item
  $C'$ is a set of representatives of those
  conjugacy classes in $G$ containing the product of  
  a diagonal matrix in $G$ with a $3$-cycle,
  and
  $\HH^2(g')$ is given in 
  Proposition~\ref{nondiagonalmatrixcontributation-nonfaithfulcohomology}.
\end{itemize}
\end{theorem}
%%%%%%%%%%%%%%%%%%%%%%%%%%%%%%%%%%%%%%%%%%%%%%%%%%%%%
\begin{proof}
By Equation (\ref{S:det1}),
$$
  \HH^2(S(V)\#G) \ = 
  \ 
   \bigoplus_{\substack{g \in {\mathcal C} \\ \codim V^g=0}}
       \HH^2(g) \ \ 
   \oplus 
   \bigoplus_{\substack{g' \in {\mathcal C} \\ \codim V^{g}=2}}
    \HH^2(g').
$$
In Proposition \ref{diagonalmatrixcontributation-nonfaithfulcohomology},
we determined the contribution to Hochschild cohomology
from those $g$ in $G$ with $\codim V^g=0$, i.e.,
from the diagonal matrices in $G$.
In Proposition \ref{nondiagonalmatrixcontributation-nonfaithfulcohomology},
we determined the contribution to Hochschild cohomology
from those $g$ in $G$ with $\codim V^g=2$, i.e.,
from those $g\in G$ conjugate to the product of a 
diagonal matrix in $G$ and a $3$-cycle.
\end{proof}

%%%%%%%%%%%%%%%%%%%%%%%%%%%%%%%%%%%%%%%%%%%%%%%%%%%%%%%%%%%%%%%%%%%%%%

\section{Graded Hecke algebras as deformations of $S(V)\# G$}
\label{gha-deformations}

In this section, we define a graded Hecke algebra associated to a
finite group $G$ and any finite dimensional representation $V$ of $G$.
We give an explicit connection between graded Hecke algebras and the
Hochschild semi-invariants defined in Section~\ref{terminologyandnotation}.
As a consequence of Theorems \ref{main-faithful-theorem} and \ref{S:paramspace}, 
if $G=G(r,p,n)$ with $r\geq 3, n\geq 4$, and $V$ is
its natural reflection representation, there are no
nontrivial graded Hecke algebras: The relevant Hochschild
cohomology, while nonzero, is not of the required form. 
This nonexistence of a graded Hecke algebra was discovered by
Ram and the first author
\cite{RamShepler}, and now we may view their result in the context
of algebraic deformation theory and Hochschild cohomology.
This negative result inspired their ad hoc construction of a
``different graded Hecke algebra'' for
$G(r,1,n)$, which coincides with an algebra
defined by Dez\'el\'ee \cite{Dezelee} in case $r=2$.
In the next section, we show that this ``different graded Hecke algebra'' is
in fact a graded Hecke algebra under our broader Definition \ref{S:GHA},
in which we allow nonfaithful actions of $G$.
This also motivates our earlier computations in Section \ref{hh2nonfaithful} of
Hochschild cohomology for a nonfaithful action of $G(r,1,n)$ on $V = \CC^n$:
Theorems \ref{main-nonfaithful-theorem} and \ref{S:paramspace} imply
existence of a generic graded Hecke algebra depending on 
many parameters, described in (\ref{S:general}).

The following definition includes
the symplectic reflection algebras defined by Chmutova
\cite{Chmutova} in case $G$ acts symplectically, i.e., when there is a
(not necessarily injective) group homomorphism $G\rightarrow {\rm Sp}(V)$.
For brevity, we omit tensor symbols in expressions in the
tensor algebra $T(V)$.

%%%%%%%%%%%%%%%%%%%%%%%%%%%%%%%%%%%
\begin{defn}\label{S:GHA}
{\em Let $G$ be a finite group with a representation $V=\CC^n$.
For each $g\in G$, choose a skew-symmetric bilinear form 
$a_g:V\times V\rightarrow \CC$. Let
$$
  A= \left(T(V)\# G\right)/(vw-wv-\sum_{g\in G}a_g(v,w)\overline{g}),
$$
where the quotient is by the ideal generated by all elements of the
given form for $v,w\in V$.
Consider $A$ to be a filtered algebra in which we assign degree 1 to
elements of $V$ and degree $0$ to elements of $G$.
We call $A$ a {\bf graded Hecke algebra} if the associated graded
algebra ${\rm {gr}} A$ is isomorphic to $S(V)\# G$.
Equivalently $A$ is isomorphic as a
vector space to $S(V)\# G$ via a choice of linear section
$S(V)\hookrightarrow T(V)$ of the canonical projection of $T(V)$ onto $S(V)$.
}\end{defn}

In case $G$ is a Coxeter group and $V$ its natural reflection
representation, the graded algebra associated to the affine Hecke algebra
(with respect to a different choice of filtration) defined by Lusztig
\cite{Lusztig3}
is a graded Hecke algebra under our definition \cite[\S3]{RamShepler}.
This is the origin of the term graded Hecke algebra.

By definition, graded Hecke algebras are parametrized by sets $\{a_g\}_{g\in G}$
satisfying suitable conditions.  
These conditions may be determined by applying \cite[Lemma~1.5]{RamShepler},
valid in this more general setting when the action of $G$ on $V$
may not be faithful:

\begin{lemma}\label{S:RS}
Let $A$ be the algebra defined by a set of
skew-symmetric bilinear forms $\{a_g\}_{g\in G}$ in Definition \ref{S:GHA}.
Then $A$ is a graded Hecke algebra if and only if
\begin{equation}\label{S:B}
  \overline{h} [v,w] (\overline{h})^{-1} = [h(v),h(w)], \ \mbox{ and}
\end{equation}

\vspace{-.7cm}

\begin{equation}\label{S:C}
  [u,[v,w]] + [v,[w,u]] + [w,[u,v]] = 0
\end{equation}
in $A$, for all $u,v,w\in V$ and $h\in G$ (where $[v,w]=vw-wv$).
\end{lemma}

Equations (\ref{S:B}) and (\ref{S:C}) are equivalent to
uniqueness of the expressions $\overline{h}wv$ and $wvu$ when rearranged
with all group elements to the right and vector space elements in
alphabetical order \cite{RamShepler}.
In case $G$ acts faithfully on $V$, this lemma was used in \cite{RamShepler}
to show directly
that the sets $\{a_g\}_{g\in G}$ corresponding to graded Hecke algebras
form a vector space of dimension $d+\dim(\Wedge^2(V))^G$, where $d$ is the
number of conjugacy classes of $g\in G$ such that $\codim(V^g)=2$ and
$\chi_g \equiv 1$ (where $\chi_g(h)=\det(h^{\perp})$ for $h\in Z(g)$).
This approach may be adapted to nonfaithful actions
to yield Corollary \ref{S:ghacor} below.
However, we take a somewhat different route in the next theorem
in order to relate graded
Hecke algebras to Hochschild cohomology, and thus to other
potential deformations of $S(V)\# G$.

\begin{defn}\label{S:defn-deg0}
{\em
Let $g\in G$ and $f_g\in\HH^2(g)$, i.e., $f_g$ 
is a Hochschild semi-invariant of $g$ (see (\ref{S:hh2g})).
By definition, we may write $f_g$ as a linear combination of elements of
the form $p\otimes y$, where $p\in S(V^g)$, $y\in\Wedge((V^g)^*)$.
We say that $f_g$ has  {\bf degree} $m$ if
the polynomial $p$ of highest degree in such a linear combination
has degree m (with respect to fixed bases of $V^g, \ \Wedge((V^g)^*)$).
Note that the homogeneous parts of $f_g$ are also Hochschild semi-invariants
of $g$.
}
\end{defn}

\begin{remark}\label{S:rem-deg0}
In the next theorem, we explain how to define a graded Hecke algebra from an element of Hochschild cohomology.  We first outline how one 
identifies a Hochschild semi-invariant with a function on $S(V)^e\ot
\Wedge^2(V)$ (recall $S(V)^e=S(V)\ot S(V)^{op}$).
Let $\mathcal C$ denote a set of representatives of the
conjugacy classes of $G$, as before.
Via (\ref{S:det1}), each element of $\HH^2(S(V)\# G)$ may be identified
with a set $\{f_g\}_{g\in {\mathcal C}}$ of Hochschild semi-invariants
(see (\ref{S:hh2g})).
Applying (\ref{S:hhsg}) as well, we regard each $f_g$ as
an element of $S(V^g)\overline{g}\otimes \Wedge ^2 (V^*)$
by making the canonical identification of the vector space
$\Wedge^{2 -\codim V^g}((V^g)^*)\otimes \Wedge^{\codim V^g}
(((V^g)^{\perp})^*)$ with a subspace of $\Wedge^{2}(V^*)$ and by
inserting the factor $\overline{g}$.
This allows us to consider $f_g$ as a $Z(g)$-invariant
function on $S(V)^e \otimes \Wedge ^2 (V)$ 
in the following way:
Suppose $f_g$ is a linear combination of elements of the form
$p\overline{g}\otimes y$ ($p\in S(V^g)$, $y\in \Wedge^2 (V^*)$).
Let $r, s\in S(V)$, $z\in \Wedge ^2 (V)$.
Then $f_g(r\otimes s\otimes z)\in S(V)\overline{g}$ is the corresponding 
linear combination of elements $y(z) r p \overline{g} s$
(see the proof of Theorem \ref{S:paramspace} 
below for details).
We further apply a transfer (trace) operator to each $f_g$ to obtain
the corresponding $G$-invariant function $f$ on $S(V)^e\ot \Wedge ^2(V)$: 
  $$f=\sum_{g\in {\mathcal C}} \ \sum_{h\in G/Z(g)} h(f_g),$$
where $(h(f_g))(r\otimes s\otimes z) = h(f_g(h^{-1}(r\otimes s\otimes z)))$.
\end{remark}

\begin{thm}\label{S:paramspace}
Let $G$ be a finite group with a representation $V = \CC^n$.
The parameter space of graded Hecke algebras for the pair $G,V$ 
is isomorphic to the space consisting of sets $\{f_g\}_{g\in G}$
of Hochschild semi-invariants
whose nonzero elements have degree~$0$.
The defining skew-symmetric bilinear forms $\{a_g\}_{g\in G}$ of the
graded Hecke algebra corresponding to 
$\{f_g\}_{g\in {{\mathcal C}}}$ are given by
$a_g(v,w)\overline{g}=f_g(1\ot 1\ot v\wedge w)$ and $a_{h^{-1}gh}(v,w)=a_g(h(v),h(w))$
for all $g\in {\mathcal C}$, $h\in G$, and $v,w\in V$.
\end{thm}

The remainder of this section is devoted to proving this theorem and its
corollary.
We need some technical lemmas
and formulas from 
\cite{CaldararuGiaquintoWitherspoon,Witherspoon1,Witherspoon2}.
We obtain the forms $a_g$ from the functions $f_g$ by finding
intermediary Hochschild two-cocycles $\mu_1$ and corresponding
deformations of $S(V)\# G$.
First, we recall the definitions of deformations and Hochschild
two-cocycles.
For more details, see \cite{Gerstenhaber} or \cite{Etingof}.

Let $t$ be an indeterminate.
If $R$ is any associative $\CC$-algebra (such as $R=S(V)\# G$),
a {\bf deformation of $R$ over $\CC[t]$} consists of the $\CC [t]$-module
$R[t]=\CC[t]\otimes R$ together with an associative product $*$ of the form
\begin{equation}\label{S:starproduct}
  a*b = ab +\mu_1(a,b)t+\mu_2(a,b)t^2+\cdots
\end{equation}
for all $a,b\in R$, where $ab$ is the product of $a$ and $b$ in $R$
and $\mu_i:R\times R\rightarrow R$ is $\CC$-bilinear (extended to be
$\CC[t]$-bilinear) for each $i$.
(In order for $a*b$ to be in $R[t]$,
this sum must in fact be finite for each pair $a,b$; 
one may also be interested in deformations over $\CC[[t]]$ or another 
extension of $\CC$, but we will not need these here.)
Associativity implies conditions on the $\mu_i$.
In particular $\mu_1$ must be a {\bf Hochschild two-cocycle}, i.e.\
\begin{equation}\label{S:2cocycle}
  \mu_1(a,bc)+a\mu_1(b,c)=\mu_1(ab,c)+\mu_1(a,b)c
\end{equation}
for all $a,b,c\in R$.
That is, $\mu_1$ is a representative of an element of $\HH^2(R)$ obtained via
the bar complex
\begin{equation}\label{S:barres}
  \cdots\stackrel{\delta_3}{\longrightarrow}
   R^{\ot 4}\stackrel{\delta_2}{\longrightarrow}
  R^{\ot 3}\stackrel{\delta_1}{\longrightarrow}
  R^{e} \stackrel{m}{\longrightarrow} R\rightarrow 0,
\end{equation}
where $\delta_i(a_0\ot a_1\ot\cdots\ot a_{i+1}) =\sum_{j=0}^i (-1)^j
    a_0\ot\cdots\ot a_ja_{j+1}\ot\cdots\ot a_{i+1}$
and $m$ is multiplication.
This is an $R^e$-free resolution of $R$ (where $R^e=R\ot R^{op}$), and
thus yields $\HH^{\DOT}(R)=\Ext^{\DOT}_{R^e}(R,R)$ upon 
taking cohomology of the cochain
complex resulting from 
application
of $\Hom_{R^e}(R, - )$.
Specifically, $\HH^i(R)=\Ker(\delta_{i+1}^*)/\Ima(\delta_i^*)$, defined
via the cochain complex
$$
0\rightarrow\Hom_{R^e}(R^e,R)\stackrel{\delta_1^*}{
  \longrightarrow} \Hom_{R^e}(R^{\ot 3},
  R) \stackrel{\delta_2^*}{
  \longrightarrow} \Hom_{R^e}(R^{\ot 4},R)\stackrel{\delta_3^*}
  {\longrightarrow} \cdots .
$$
We identify $\Hom_{R^e}(R^{\ot 4},R)\cong \Hom_{\CC}(R^{\ot 2},R)$, and
a straightforward calculation yields Equation (\ref{S:2cocycle}) as
the defining relation for elements of $\Ker(\delta_3^*)$.

In order to prove Theorem \ref{S:paramspace}, we explain how to
obtain the defining skew-symmetric forms $\{a_g\}_{g\in G}$ 
of a graded Hecke algebra from a set of
Hochschild semi-invariants $\{f_g\}_{g\in {\mathcal C}}$ whose
nonzero elements have degree $0$.
The following comparison of the bar complex (\ref{S:barres}) to a Koszul
complex allows us to write down corresponding Hochschild
2-cocycles $\mu_1$ explicitly.

There is a chain map from the bar complex (\ref{S:barres}) for 
$R=S(V)$ to the Koszul complex
$P_{\DOT} = K(\{v_i\ot 1 - 1\ot v_i\}_{i=1}^n)$, where $v_1,\ldots,v_n$ is a
basis of $V$ (see \cite[\S4.5]{Weibel} for details on Koszul complexes):
$$
\begin{array}{ccccccccccc}
\cdots \! &\rightarrow & S(V)^{\ot 4} & \stackrel{\delta_2}{\longrightarrow}
  & S(V)^{\ot 3} & \stackrel{\delta_1}{\longrightarrow} & S(V)^e &
  \stackrel{m}{\longrightarrow} & S(V) & \rightarrow & \! 0\\
 & & \hspace{.1in}\downarrow \psi_2 & & \hspace{.1in}\downarrow \psi_1 & &
   \parallel &&\parallel &&\\
 \cdots \! &\! \rightarrow \! & S(V)^e\ot \Wedge^2(V) \! & \! 
  \stackrel{d_2}{\longrightarrow}\!
  & \!S(V)^e\ot \Wedge^1(V) \! &\! \stackrel{d_1}{\longrightarrow} 
  \! & \! S(V)^e &\stackrel
  {m}{\longrightarrow}\! & \! S(V) \! &\! \rightarrow & \! 0 
\end{array}
$$
The differentials $d_1,d_2$ are  $S(V)^e$-homomorphisms given in our notation
by 
\begin{eqnarray*}
d_1(1\ot 1\ot v_i)&=&v_i\ot 1-1\ot v_i , \\
d_2(1\ot 1\ot v_i\wedge v_j)&=& (v_i\ot 1 - 1\ot v_i)\ot v_j
- (v_j\ot 1 -1\ot v_j)\ot v_i,
\end{eqnarray*}
for $1\leq i<j\leq n$.
The vertical maps $\psi_1$ and $\psi_2$ may be given by the formulas 
(see \cite[(4.9) and (4.10)]{Witherspoon1}, where the notation is somewhat
different):
\begin{equation}\label{eqn:psi1}
 \psi_1(1\ot v_1^{k_1}\cdots v_n^{k_n}\ot 1) = \sum_{i=1}^n \sum_{a=1}
   ^{k_i} v_i^{k_i-a}v_{i+1}^{k_{i+1}}\cdots v_n^{k_n}\ot 
  v_1^{k_1}\cdots v_{i-1}^{k_{i-1}} v_i^{a-1} \ot v_i,
\end{equation}
\begin{equation}\label{S:psi2}
\psi_2(1\ot v_1^{k_1}\cdots v_n^{k_n}\ot v_1^{m_1}\cdots v_n^{m_n}\ot 1)=
\hspace{4in}
\end{equation}
$$\sum_{1\leq i<j\leq n}\sum_{b=1}^{m_j}\sum_{a=1}^{k_i} 
  v_i^{k_i-a} v_{i+1}^{k_{i+1}}\cdots v_{j-1}^{k_{j-1}}v_j^{k_j+m_j-b}v_{j+1}
  ^{k_{j+1}+m_{j+1}}\cdots v_n^{k_n+m_n}\ot 
$$

\vspace{-.15in}

$$\hspace{1.5in}v_1^{k_1+m_1}\cdots v_{i-1}^
 {k_{i-1}+m_{i-1}} v_i^{m_i+a-1} v_{i+1}^{m_{i+1}}\cdots v_{j-1}^{m_{j-1}}
  v_j^{b-1} \ot v_i\wedge v_j.
$$
To obtain an explicit Hochschild two-cocycle $\mu_1$ from a set
$\{f_g\}_{g\in {{\mathcal C}}}$ of Hochschild semi-invariants, we apply
the following proposition. This
proposition appears as part of
\cite[Thm.\ 5.4]{CaldararuGiaquintoWitherspoon},
where it is stated for arbitrary degree and arbitrary resolution $P_{\DOT}$. It is
valid also when $G$ acts nonfaithfully on $V$.

%%%%%%%%%%%%%%%%%%%%%%%%%%%%%%%%%%%%%%%%%%
\begin{prop}\label{S:CGW}
Let $R=S(V)\# G$.
Let $f:S(V)^e\ot \Wedge^2(V)\rightarrow R$ be a function
(on the degree 2 term of the above Koszul complex) 
representing an element of $\HH^2(S(V),R)^G$.
Under the isomorphism $\HH^2(S(V),R)^G\cong \HH^2(R)$, 
$f$ is mapped to the function 
$\mu_1\in \Hom_{\CC}(R^{\ot 2},R)\cong
\Hom_{R^e}(R^{\ot 4}, R)$ whose action
on the degree 2 term of the bar complex (\ref{S:barres}) is given by
$$
  \mu_1 (r\overline{g}\ot s\overline{h})=((f\circ \psi_2)
 (1\ot r\ot g(s)\ot 1))\overline{gh},
$$
for all $r,s\in S(V)$ and $g,h\in G$.
\end{prop}
%%%%%%%%%%%%%%%%%%%%%%%%%%%%%%%%%%%%%%%%%%

Our final tool for proving Theorem \ref{S:paramspace} is
\cite[Thm.\ 3.2]{Witherspoon2}, valid when
$G$ acts nonfaithfully on $V$; we also record it here for convenience.
Consider $S(V)\# G$ to be a graded algebra where elements of $V$
have degree $1$ and elements of $G$ have degree~$0$.
In \cite{Witherspoon2}, a graded Hecke algebra is defined
over $\CC[t]$: It is a quotient of $T(V)\# G [t]$, by the ideal 
generated by all $vw-wv-\sum_{g\in G} a_g(v,w)t\overline{g}$,
whose associated graded $\CC[t]$-algebra is isomorphic to
$S(V)\# G[t]$.
To obtain our Definition \ref{S:GHA}, simply substitute any
nonzero complex number for $t$.

\begin{prop}\label{S:W}
Up to isomorphism, the (nontrivial) graded Hecke algebras over $\CC[t]$ are precisely the
deformations of $S(V)\# G$ over $\CC[t]$ for which $\deg \mu_i = -2i$
($i\geq 1$).
\end{prop}

\begin{proof}[Proof of Theorem \ref{S:paramspace}]
Suppose $\{f_g\}_{g\in {\mathcal C}}$ is a set of Hochschild semi-invariants
and that each nonzero $f_g$ has degree $0$.
As in Remark \ref{S:rem-deg0}, 
let $f$ be the associated function on $S(V)^e\ot \Wedge^2(V)$.
Let $\mu_1$ be the corresponding Hochschild 2-cocycle of $S(V)\# G$, given
by Proposition \ref{S:CGW}.
Let $v,w\in V$.
By Remark \ref{S:rem-deg0}, Proposition \ref{S:CGW}, and formula
(\ref{S:psi2}), $\mu_1(v\ot w)=(f\circ \psi_2)(1\ot v\ot w\ot 1)$
is an element of the subalgebra
$\CC G$ of $S(V)\# G$.
Similarly, $\mu_1(w\ot v)$ is an element of $\CC G$.
Define scalars $a_g(v,w)$ by setting
\begin{equation}\label{S:A}
  \mu_1(v\ot w)-\mu_1(w\ot v) = \sum_{g\in G} a_g(v,w)\overline{g}.
\end{equation}
As $\mu_1$ is bilinear, each resulting function $a_g: V\times V\rightarrow \CC$
is bilinear.
By definition, $a_g$ is skew-symmetric for each $g\in G$.
We claim that the set $\{a_g\}_{g\in G}$ defines a graded Hecke algebra.
Let $A$ be the corresponding quotient given in Definition~\ref{S:GHA}.
By Lemma \ref{S:RS}, it suffices to verify (\ref{S:B}) and (\ref{S:C}) for $A$.

By Definition \ref{S:GHA} and Equation (\ref{S:A}), 
$[v,w]=\mu_1(v\ot w)-\mu_1(w\ot v)$ lies in $A$,
for all $v,w\in V$.
Abuse notation and write $\psi_2(v\ot w)=\psi_2(1\ot v\ot w\ot 1)$.
Apply definition (\ref{S:psi2}) of $\psi_2$ 
to our chosen basis for $V$ to see that
$\psi_2(v_i\ot v_j)=1\ot 1\ot v_i\wedge v_j$ if $i<j$ and $0$ otherwise.
After expressing all vectors in terms of this basis, we find
$$
  h\left(\psi_2(v\ot w)-\psi_2(w\ot v)\right) =\psi_2(h(v)\ot h(w))-
  \psi_2(h(w)\ot h(v))
$$
for all $v,w\in V$, $h\in G$. 
As $f$ is also $G$-invariant, it follows that 
$$\mu_1(v\ot w)-\mu_1(w\ot v)=
(f\circ \psi_2)(v\ot w)-(f\circ\psi_2)(w\ot v)
$$ 
is $G$-invariant 
as a function on $v,w\in V$, and thus (\ref{S:B}) holds in $A$.

Next we claim that the Jacobi identity (\ref{S:C}) 
is a direct consequence of
the Hochschild 2-cocycle condition (\ref{S:2cocycle}).
The left side of (\ref{S:C}) may be rewritten by replacing the innermost
bracket $[v,w]$ in the first term by $\mu_1(v\ot w)-\mu_1(w\ot v)$, and
similarly for each of the other two terms.  We obtain
$$
  u\mu_1(v\ot w)-u\mu_1(w\ot v)-\mu_1(v\ot w)u+\mu_1(w\ot v)u+v\mu_1(w\ot
  u)-v\mu_1(u \ot w)
$$
$$
  \hspace{.5cm} -\mu_1(w\ot u)v+\mu_1(u\ot w)v +w\mu_1(u\ot v)-w\mu_1(v\ot u)
   -\mu_1(u\ot v)w+\mu_1(v\ot u)w.
$$
Substitutions from the six permutations of (\ref{S:2cocycle}) in which
$\{a,b,c\}=\{u,v,w\}$ yield
$$
  \mu_1(uv\ot w)-\mu_1(u\ot vw)+\mu_1(u\ot wv)-\mu_1(uw\ot v)+\mu_1(vw\ot u)
  -\mu_1(v\ot wu)
$$
$$
  \hspace{.5cm} +\mu_1(w\ot vu)-\mu_1(wv\ot u)+\mu_1(v\ot uw)-\mu_1(vu\ot w)
  +\mu_1(wu\ot v)-\mu_1(w\ot uv).
$$
But this expression is zero as $S(V)$ is commutative.
We have thus shown that (\ref{S:C}) holds in $A$.
This concludes the proof that a set of
Hochschild semi-invariants whose nonzero elements have
degree $0$  gives rise to a graded Hecke algebra.

We now argue that every graded Hecke algebra arises in this way.
Note that if a set $\{a_g\}_{g\in G}$ of skew-symmetric bilinear forms
defines a graded Hecke algebra, then so does $\{\alpha  a_g\}_{g\in G}$
for any fixed scalar $\alpha$.
We may therefore consider the related graded Hecke algebra 
over $\CC[t]$, where the indeterminate
$t$ takes the place of the arbitrary scalar $\alpha$. 
By Proposition \ref{S:W}, a (nontrivial) graded Hecke algebra over $\CC[t]$
is a deformation of $S(V)\# G$ such that the associated
Hochschild 2-cocycle $\mu_1$ satisfies $\deg \mu_1 =-2$, and more generally
$\deg(\mu_i)=-2i$ ($i\geq 1$). 
Equation (\ref{S:starproduct}) and Definition \ref{S:GHA} (with $a_g$
replaced by $ta_g$)
force relationship (\ref{S:A}) between the defining forms $ta_g$  and the cocycle $\mu_1$.
As $\mu_1$ is a Hochschild 2-cocycle,
Equation (\ref{S:det1}) implies existence of a set of
Hochschild semi-invariants $\{f_g\}_{g\in {\mathcal C}}$ and a corresponding 
function $f$ (see Remark \ref{S:rem-deg0}) such that 
$\psi_2^*(f)$ is cohomologous to $\mu_1$. 
In particular, there exists an $S(V)^e$-homomorphism
$\beta :S(V)^{\ot 3} \rightarrow S(V)\# G$ such that 
$f\circ \psi_2 = \mu_1+ \beta\circ\delta_2$.
Identify $\beta$ with a $\CC$-linear function from $S(V)$ to
$S(V)\# G$, and note that
$$\beta(\delta_2(1\ot v\ot w\ot 1))
=v\beta(w) -\beta(vw)+\beta(v)w$$
for all $v,w\in V$.
Thus for all $i<j$ we have
$$
  (f\circ \psi_2)(v_i\ot v_j)=\mu_1(v_i\ot v_j)-\mu_1(v_j\ot v_i)+v_i\beta(v_j)
  +\beta(v_i)v_j -v_j\beta(v_i)-\beta(v_j)v_i
$$
(since $\psi_2(v_j\ot v_i)=0$ when $i<j$).
Thus $(f\circ \psi_2)(v_i\ot v_j)$ is a sum of the element $\mu_1(v_i\ot v_j)-
\mu_1(v_j\ot v_i)$ in $\CC G$ and the element 
$v_i\beta(v_j)+\beta(v_i)v_j -v_j\beta(v_i)-\beta(v_j)v_i$
in the ideal $(V)$ of
$S(V)\# G$. 
Now $(f\circ \psi_2)(v_i\ot v_j)=f(1\ot 1\ot v_i\wedge v_j)$, and
the function $f$ is discussed in Remark \ref{S:rem-deg0};
in particular it is determined by these values since it is an 
$S(V)^e$-homomorphism.
Since $S(V)\# G \cong \CC G \oplus S(V)$ as a vector space
and the homogeneous parts of $f$ are also $G$-invariant functions on
$S(V)^e\ot \Wedge ^2(V)$, we have $f=f'+f''$ where $(f'\circ \psi_2)(v_i\ot v_j)
=\mu_1(v_i\ot v_j)-\mu_1(v_j\ot v_i)$ and 
$f''\circ \psi_2 = \beta\circ\delta_2$ is a coboundary.
By Remark \ref{S:rem-deg0}, $f'$ itself corresponds to a set of Hochschild semi-invariants,
$\{f_g'\}_{g\in {\mathcal C}}$. As $\deg \mu_1=-2$ and $\deg \psi_2 =-2$,
each $f_g'$ is either 0 or has degree $0$.
Thus $\{f_g'\}_{g\in {\mathcal C}}$ is the desired set of Hochschild semi-invariants.
(Alternatively, we could simply prove that our map from sets
$\{f_g\}_{g\in {\mathcal C}}$ of Hochschild semi-invariants (whose nonzero
elements have degree $0$) to sets $\{a_g\}_{g\in G}$ of skew-symmetric
forms defining graded Hecke algebras is {\em surjective} by classifying
the latter directly from Lemma \ref{S:RS}, cf.\ \cite[Thm.\ 1.9]{RamShepler}.) 

The formula for the $a_g$ in terms of the $f_g$ follows from our calculations
above:
If $i<j$, then
\begin{eqnarray*}
  \sum_{g\in G}a_g(v_i,v_j)\overline{g} \ = \ [v_i,v_j]
   &=& \mu_1(v_i\ot v_j)-\mu_1(v_j\ot v_i)\\
  &=& (f\circ \psi_2)(v_i\ot v_j) \ = \ f(1\ot 1\ot v_i\wedge v_j).
\end{eqnarray*}
Now examine the coefficient of each $\overline{g}$ ($g\in G$) separately,
and apply (\ref{S:B}).
The resulting formula
extends from pairs $v_i,v_j$ of basis elements to all $v,w\in V$
by linearity.
\end{proof}

\begin{remark}\label{S:Jacobi}
We point out an interesting consequence of the Hochschild
2-cocycle condition (\ref{S:2cocycle}) in this context:
Let $\mu_1$ and $A$ be as in the first paragraph
of the proof of Theorem \ref{S:paramspace}.
We showed that the Jacobi identity (\ref{S:C}) holds in $A$ 
as a direct consequence of the Hochschild 2-cocycle condition
(\ref{S:2cocycle}) for $\mu_1$.
\end{remark}

%%%%%%%%%%%%%%%%%%%%%%%%%%%%%%%%%%%%%%%%%%%%%%%%%%%%%%%%%%%%%%%%%%%%%%

The following corollary characterizes graded Hecke algebras without
reference to Hochschild cohomology.
Recall that $\chi_g(h)=\det(h^{\perp})$ for $h\in Z(g)$, 
where $h^{\perp}=h|_{(V^g)^{\perp}}$,
and $\mathcal C$ denotes a set of representatives of the conjugacy classes of $G$.

\begin{cor}\label{S:ghacor}
Let $d$ be the number of conjugacy classes of $g\in G$ such that
$\codim V^g =2$ and $\chi_g \equiv 1$.
\begin{itemize}
\item[(i)] The sets $\{a_g\}_{g\in G}$ corresponding to graded Hecke algebras
form a vector space of dimension 
  $$d \ \ +\sum_{\substack{g\in {\mathcal {C}} \\  V^g = V}} 
      \dim (\Wedge^2(V))^{Z(g)}.$$
\item[(ii)] All graded Hecke algebras arise in the following way:
For each $g\in {\mathcal C}$ satisfying $\codim V^g =2$ and $\chi_g \equiv 1$, 
define a skew-symmetric bilinear form $a_g:V\times V\rightarrow\CC$
by setting $a_g(w_1,w_2)$ equal to an arbitrary scalar
for a fixed basis $\{w_1,w_2\}$ of $(V^g)^{\perp}$, $a_g(V^g,V)=0$, and
$a_{h^{-1}gh}(v,w)=a_g(h(v),h(w))$ for all $h\in G, v,w\in V$.
For each $g\in {\mathcal C}$ satisfying $V^g =V$, define $a_g$ 
by 
any choice of $Z(g)$-invariant linear functional on $\Wedge^2(V)$,
and $a_{h^{-1}gh}(v,w)=a_g(h(v),h(w))$ for all $h\in G$, $v,w\in V$.
\end{itemize}
\end{cor}

\begin{proof}
By (\ref{S:det1}), if $\HH^2(g)$ is nonzero, then $\codim V^g\in \{0,2\}$.
We apply Theorem~\ref{S:paramspace} in each case:
If $\codim V^g =0$ (i.e.\ $V^g=V$), 
the Hochschild semi-invariants of $g$ that are either 0 or
of degree $0$ form a vector space of dimension 
$\dim (\Wedge^2(V))^{Z(g)}$ (see (\ref{S:hh2g})). 
Theorem \ref{S:paramspace}, Remark \ref{S:rem-deg0}, and 
(\ref{S:B}) then give the indicated
form of the corresponding functions $a_g$.
If $\codim V^g =2$, a nonzero scalar is a Hochschild semi-invariant
for $g$ if and only if $\chi_g\equiv 1$ (by definition).
Again by Theorem \ref{S:paramspace}, Remark \ref{S:rem-deg0}, and (\ref{S:B}),
the values of the corresponding forms $a_g$ are precisely those given in (ii).
(Alternatively the corollary may be proven directly, without using
Hochschild cohomology, cf.\ \cite[Thm.\ 1.9]{RamShepler}). 
\end{proof}

\section{Graded Hecke algebras for $G(r,1,n)$}
\label{ghaG(r,1,n)}

Assume $n\geq 3$. Let $G= G(r,1,n)$ act on $V=\CC^n$ nonfaithfully as the symmetric group via $\rho$
defined in (\ref{S:rep}).
In this section, we find all graded Hecke algebras corresponding to this
action. Similar results hold for $G(r,p,n)$.
We prove that in a special case, these algebras are
isomorphic to algebras that appeared in
\cite{RamShepler}, and in \cite{Dezelee} when $r=2$.

By Theorems \ref{main-nonfaithful-theorem} and \ref{S:paramspace},
the parameter space of graded Hecke algebras for $G,V$
has dimension equal to the number of conjugacy classes 
of elements $g$ in $G(r,1,n)$ for which $g$ is 
the product of a diagonal matrix with a 3-cycle.
(By Proposition \ref{diagonalmatrixcontributation-nonfaithfulcohomology}
and Remark \ref{diagonalmatrixcontributation-nonfaithfulcohomology-remark},
the diagonal matrices themselves do not have Hochschild semi-invariants
of degree $0$.)
We refine the notation introduced in the text preceding
Proposition \ref{nondiagonalmatrixcontributation-nonfaithfulcohomology}
for these elements $g$: Choose representatives $g=g_{b,\underline{c},
\underline{m}}$, one for each conjugacy class, where
$$
  g_{b,\underline{c},\underline{m}}=\xi_3^b(1,2,3)\xi^{c_1}I_{m_1}\oplus
\cdots\oplus\xi^{c_k}I_{m_k},
$$
$\underline{c}=(c_1,\ldots,c_k)$, $\underline{m}=(m_1,\ldots,m_k)$, and
$0\leq b , c_1,\ldots,c_k\leq r-1$ ($m_1,\ldots,m_k$ are the multiplicities
of the diagonal entries other than the first three).
By Corollary \ref{S:ghacor}(ii), all graded Hecke algebras arise as follows.
For each representative $g_{b,\underline{c}, \underline{m}}$, choose a scalar 
$\alpha_{b,\underline{c}, \underline{m}}$.
The corresponding graded Hecke algebra is
\begin{equation}\label{S:general}
  A=\left(T(V)\# G\right)/(vw-wv-\sum_{g\in G(r,1,n)} a_g(v,w)
  \overline{g}),
\end{equation}
where the skew-symmetric bilinear forms $a_g$ ($g\in G$) are given by:
\begin{itemize}
\item $a_g\equiv 0$ if
$g$ is not conjugate to one of the representatives $g_{b,\underline{c},
\underline{m}}$; 
\item $a_{g_{b,\underline{c}, \underline{m}}}
  (v_1-v_2, \ v_2-v_3)=\alpha_{b,\underline{c}, \underline{m}} \ $
and 
$ \  a_{g_{b,\underline{c}, \underline{m}}}
(V^{(1, 2, 3)},V)=0$;
\item $a_{h^{-1}g_{b,\underline{c}, \underline{m}}h}
  (v,w) = 
 a_{g_{b,\underline{c}, \underline{m}}}
  (h(v),h(w)) \ $ for all $h\in G$.
\end{itemize}

We are particularly interested in the choice $\alpha_{0,0,n-3}=1$ ($k=1$) and
all other $\alpha_{b,\underline{c},\underline{m}}=0$. Denote the
resulting graded Hecke algebra by $A_{r,1,n}$:

\begin{defn}\label{S:gha}
 $ \ \displaystyle{ A_{r,1,n}= \left(T(V)\# G(r,1,n)\right) /
  (vw-wv-\sum_{g\in G(r,1,n)} a_g(v,w) \overline{g}),}$

where the skew-symmetric bilinear forms $a_g$ ($g\in G$) are given by:
\begin{itemize}
\item $a_g\equiv 0$ if $g$ is not conjugate to $(1, 2, 3)$;
\item $a_{(1, 2, 3)}(v_1-v_2 , v_2-v_3)=1$  
and $a_{(1, 2, 3)}(V^{(1, 2, 3)} , V)=0$;
\item  $a_{h^{-1}(1, 2, 3)h}(v,w)=a_{(1, 2, 3)}(h(v),h(w))$ for all
$h\in G$.
\end{itemize}
\end{defn}

We claim
that the defining relations $vw-wv-\sum_{g\in G}a_g(v,w)
\overline{g}$ of $A_{r,1,n}$ may be replaced by
%%%%%%%%%%%%%%%%%%%%%%%%%%%%%%%%%%%%%%%%%%%%%%%%%%%%%%%%%%%%%%%%%%%%%%%
\begin{equation}\label{S:89}
  v_mv_k-v_kv_m=\frac{1}{3} 
  \sum_{\substack{1\leq i\leq n \\ i\neq m,k}}
  \ \sum_{a , b =0}^{r-1}
  \overline{\xi}_m^{a}\overline{\xi}_k^b
  (\overline{\xi}_i)^{-a -b} ((\overline{m, k, i})-
  (\overline{m, i, k}))
\end{equation}
for all $1\leq m<k\leq n$.
To see this, first note that the group elements in this sum are precisely
those that both are conjugate to $(1,2,3)$ and act nontrivially on $v_m,v_k$.
For each $i,a,b$ let $h$ be an element of $\Sym_n$ such that
$h^{-1}(1,2,3)h=\xi_m^{a}\xi_k^b\xi_i^{-a -b}(m,k,i)$ 
(or $h^{-1}(1,2,3)h=\xi_m^{a}\xi_k^b\xi_i^{-a -b}(m,i,k)$).
Then apply  Definition \ref{S:gha}, writing 
$v_1,v_2,v_3$ in terms of the vector
space decomposition $V\cong (V^{(1,2,3)})^{\perp}\oplus V^{(1,2,3)}$.

We will show that the graded Hecke algebra $A_{r,1,n}$
is precisely the algebra defined by Ram and the first author 
\cite[(5.1)]{RamShepler}
as a substitute for a (nonexistent) graded Hecke algebra corresponding
to the natural reflection representation of $G(r,1,n)$.
Their algebra agrees with one defined by Dez\'el\'ee \cite{Dezelee} in case $r=2$.
Our results in Section \ref{gha-deformations} will then imply
that these algebras of \cite{Dezelee,RamShepler}
arise from deformations of $S(V)\# G(r,1,n)$ over $\CC [t]$.
They were originally given in a different form, analogous to
Lusztig's definition of a graded Hecke algebra.
The proof that Lusztig's definition is a special case of Drinfeld's 
definition \cite[Thm.\ 3.5]{RamShepler} suggested to us how to proceed.

\begin{defn}\label{S:hstar}{\em
As in \cite[(5.1)]{RamShepler},
let $H^*_{r,1,n}$ be the algebra generated by the basis $v_1,\ldots,v_n$
of $V$ and all $\overline{g}$ ($g\in G(r,1,n)$) such that
$\CC G(r,1,n)$ and $S(V)$ are subalgebras and the following additional 
relations hold:
\begin{eqnarray}
  \overline{\xi}_i v_k &=& v_k \overline{\xi}_i  \hspace{1.5in} 
   (1\leq i,k\leq n), \label{S:reln1}\\
  \overline{s}_i v_k& =& v_k \overline{s}_i \hspace{1.51in} (k\neq i,i+1),
  \label{S:reln2}\\
  \overline{s}_i v_{i+1} &=& v_i\overline{s}_i + \sum_{a = 0}^{r-1}\
   (\overline{\xi}_i)^{\,a} (\overline{\xi}_{i+1})^{-a} \hspace{.385in}
   (1\leq i\leq n-1),\label{S:reln3}
\end{eqnarray}
where $s_i$ is the simple reflection $(i,i+1)$.
}\end{defn}

A useful variant of relation (\ref{S:reln3}) is $\overline{s}_iv_i=v_{i+1}
\overline{s}_i - \sum_{a =0}^{r-1}\ (\overline{\xi}_i)^{\,a}\, 
(\overline{\xi}_{i+1})^{-a}$.
In Dez\'el\'ee's version \cite{Dezelee}, one may choose any nonzero scalar as the
coefficient of the sum in relation (\ref{S:reln3})
above, and the proof of the following theorem may be modified accordingly.

\begin{thm} \label{S:iso}
Let $A_{r,1,n}$ be the graded Hecke algebra of 
Definition \ref{S:gha}.
There is an algebra isomorphism
$$
  A_{r,1,n}\cong H^*_{r,1,n}.
$$
\end{thm}

The remainder of this section is devoted to the proof of this theorem.
We first collect some additional relations in $H^*_{r,1,n}$
that are consequences of (\ref{S:reln2}) and (\ref{S:reln3}).

\begin{lemma}\label{S:reln4}
Let $1\leq j,k,m \leq n$ with $j<k$.
In $H^*_{r,1,n}$, the element $(\overline{j, k}) v_m$ is equal to
\begin{equation*}
\left\{
\begin{array}{ll}
  v_m(\overline{j, k}), & \mbox{if } m<j \mbox{ or } k<m\\
  \displaystyle{v_m(\overline{j, k})+
  \sum_{a =0}^{r-1} (\overline{\xi}_m)^{\,a}
   (\overline{\xi}_k)^{-a}\, (\overline{j, m, k}) -\sum_{a =0}^{r-1}
   (\overline{\xi}_j)^{\,a} (\overline{\xi}_m)^{-a}\, 
   (\overline{j, k, m})},
   & \mbox{if } j<m<k\\
  \displaystyle{v_k(\overline{j, k}) 
  - \sum_{j<i<k} \sum_{a =0}^{r-1} (\overline{\xi}_i)
    ^{\,a}(\overline{\xi}_k)^{-a}\, (\overline{j, i, k})
    -\sum_{a =0}^{r-1} (\overline{\xi}_j)^{\,a}\,
   (\overline{\xi}_k)^{-a}} , & \mbox{if } m=j\\
  \displaystyle{v_j(\overline{j, k}) 
   + \sum_{j<i<k}\sum_{a =0}^{r-1} (\overline{\xi}_i)
  ^{\,a} (\overline{\xi}_j)^{-a}\, (\overline{k, i, j})
  +\sum_{a =0}^{r-1} (\overline{\xi}_j)^{\,a}(\overline{\xi}_k)^{-a}}, &
   \mbox{if } m=k.
\end{array}\right.
\end{equation*}
\end{lemma}

\begin{proof}
Write $(\overline{j, k})$ as a product of simple reflections:
$$
  (\overline{j, k}) = \overline{s}_{k-1}\overline{s}_{k-2}\cdots
   \overline{s}_{j+1}\overline{s}_j \overline{s}_{j+1}\cdots
   \overline{s}_{k-2}\overline{s}_{k-1}.
$$
If $m<j$ or $k<m$, then $(\overline{j, k}) v_m = v_m (\overline{j, k})$ by
(\ref{S:reln2}).
If $j<m<k$, then application of (\ref{S:reln3}) and its variation
(given in the text following (\ref{S:reln3}))
yields the stated relation.

If $m=j$, we use induction on $k-j$. If $k-j=1$ then the desired
relation is the one given in the text following (\ref{S:reln3}).
Assume the
relation given in Lemma \ref{S:reln4} holds when $k$ is replaced by
$k-1$ and $m=j<k-1$. Then
 $ (\overline{j, k}) v_j = \overline{s}_{k-1} (\overline{j, k\! -\! 1}) 
  \overline{s}_{k-1} v_j
  = \overline{s}_{k-1} (\overline{j, k\! -\! 1}) v_j \overline{s}_{k-1}$,
and application of the induction hypothesis and the variation of
(\ref{S:reln3}) yields the relation stated in the theorem.
The final relation is proved by a similar induction.
\end{proof}

We will make a change of generators for $H^*_{r,1,n}$. 
For each $k$ ($1\leq k\leq n$), let
\begin{equation}\label{S:wktilde}
\widetilde{v}_k = v_k + \frac{1}{2}\sum_{\substack{1\leq j\leq n \\
   j\neq k}}\sum_{a =0}^{r-1}\ (-1)^{\delta_{j<k}}\
  (\overline{\xi}_k)^{\, a}\ (\overline{\xi}_j)^{-a}\ (\overline{k, j}),
\end{equation}
where $\delta_{j<k} =1$ if $j<k$ and $0$ otherwise.

\begin{lemma}
The algebra $H^*_{r,1,n}$ is generated by $\widetilde{v}_k$ ($1\leq k\leq n$)
and $\overline{g}$ ($g\in G(r,1,n)$) and is defined by the relations of
$G(r,1,n)$ together with
\begin{eqnarray}
 \overline{\xi}_i\widetilde{v}_k & = & \widetilde{v}_k \overline{\xi}_i \ \ \ 
   \mbox{ for all }1\leq i,k\leq n,\label{S:trivial1}\\
 \overline{s}_i\widetilde{v}_k &=& \widetilde{v}_k \overline{s}_i \ \ \ 
   \mbox{ if } k\not\in \{i,i+1\},\label{S:trivial2}\\
 \overline{s}_i\widetilde{v}_i & =& \widetilde{v}_{i+1} \overline{s}_i \ \ \
  \mbox{ for }1\leq i\leq n-1,\label{S:action}\\
 \quad\quad\widetilde{v}_m\widetilde{v}_k - \widetilde{v}_k
  \widetilde{v}_m &=& 
\frac{1}{4}\sum_{\substack{1\leq i\leq n \\ i\neq m,k
 \rule[0ex]{0ex}{1.5ex}%strut
}}
   \ \sum_{a , b=0}^{r-1}(\overline{\xi}_m)^{a}
   (\overline{\xi}_k)^b (\overline{\xi}_i)^{-a -b}
   ((\overline{m, k, i})-(\overline{m, i, k})),\label{S:bracket}
\end{eqnarray} 
where the last relation holds for $1\leq m<k\leq n$.
\end{lemma} 

\begin{proof}
In the definition of $\widetilde{v}_k$, 
both summands
$v_k$ and $\sum \sum (-1)^{\delta_{j<k}} (\overline{\xi}_k)^{a}(\overline{
\xi}_j)^{-a} (\overline{k, j})$ are invariant
under conjugation by $\overline{\xi}_i$ ($1\leq i\leq n$) and by
$\overline{s}_i$ ($k\neq i,i+1$), so (\ref{S:trivial1}) and 
(\ref{S:trivial2}) hold.

We check (\ref{S:action}), using (\ref{S:reln3}):
\begin{eqnarray*}
\overline{s}_i\widetilde{v}_i\overline{s}_i &=& \overline{s}_i v_i
   \overline{s}_i + \frac{1}{2}\sum_{\substack{1\leq j\leq n \\
   j\neq i
   \rule[0ex]{0ex}{1.5ex}%strut
  }}\ \sum_{a =0}^{r-1} (-1)^{\delta_{j<i}}\, \overline{s}_i
   (\overline{\xi}_i)^{a}(\overline{\xi}_j)^{-a} (\overline{i, j})
  \overline{s}_i\\
 &=& v_{i+1} -\sum_{a =0}^{r-1} \overline{s}_i (\overline{\xi}_i)^{a}
  (\overline{\xi}_{i+1})^{-a}  +\frac{1}{2}\sum_{\substack{1\leq j\leq n \\
    j\neq i\! +\! 1
    \rule[0ex]{0ex}{1.5ex}%strut
   }}\ \sum_{a =0}^{r-1} (-1)^{\delta_{j<i}} 
  (\overline{\xi}_{i+1})^{a} (\overline{\xi}_j)^{-a} (\overline{i\! +\! 
    1, j})\\
 &=& v_{i+1} +\frac{1}{2}\sum_{\substack{1\leq j\leq n \\ j\neq i+1
  \rule[0ex]{0ex}{1.5ex}%strut
  }}\
   \sum_{a =0}^{r-1} (-1)^{\delta_{j< (i\! +\! 1)}}
  (\overline{\xi}_{i+1})^{a}
  (\overline{\xi}_j)^{-a} (\overline{i\! +\! 1, j})\\
  &=& \widetilde{v}_{i+1}.
\end{eqnarray*}

Finally, we check (\ref{S:bracket}). Assume $m<k$.
We compute the bracket $\widetilde{v}_m\widetilde{v}_k -\widetilde{v}_k
\widetilde{v}_m$, first substituting from (\ref{S:wktilde}) and then
applying the relations from Lemma \ref{S:reln4}.
The result is a cancelation
of all terms having factors of $v_m$ or $v_k$:
$$
(v_m+\frac{1}{2}\sum_{\substack{1\leq i\leq n \\ i\neq m
\rule[0ex]{0ex}{1.5ex}%strut
  }}
  \sum_{b=0}^{r-1} (-1)^{\delta_{i<m}}(\overline{\xi}_m)^b(\overline{\xi}_i)^{-b}
  (\overline{m, i}))(v_k+\frac{1}{2}\sum_{\substack{1\leq j\leq n \\ j\neq k
   \rule[0ex]{0ex}{1.5ex}%strut
    }}
   \sum_{a =0}^{r-1}(-1)^{\delta_{j<k}}(\overline{\xi}_k)^{a}
   (\overline{\xi}_j)^{-a}(\overline{k, j}))
$$
$$
-(v_k+\frac{1}{2}\sum_{\substack{1\leq j\leq n \\ j\neq k
   \rule[0ex]{0ex}{1.5ex}%strut
   }}
   \sum_{a =0}^{r-1}(-1)^{\delta_{j<k}}(\overline{\xi}_k)^{a}
   (\overline{\xi}_j)^{-a}(\overline{k, j}))(v_m+\frac{1}{2}
   \sum_{\substack{1\leq i\leq n \\ i\neq m
   \rule[0ex]{0ex}{1.5ex}%strut
   }}
   \sum_{b=0}^{r-1} (-1)^{\delta_{i<m}}(\overline{\xi}_m)^b
  (\overline{\xi}_i)^{-b}(\overline{m, i}))
$$

$$\hspace{-.3cm}
= \frac{1}{2}\sum_{\substack{1\leq j\leq n \\ j\neq k
   \rule[0ex]{0ex}{1.5ex}%strut 
   }}\sum_{a =0}^{r-1}
  (-1)^{\delta_{j<k}}v_m (\overline{\xi}_k)^{a}(\overline{\xi}_j)^{-a}
  (\overline{k, j}) +\frac{1}{2} \sum_{\substack{1\leq i\leq n \\ i\neq m
  \rule[0ex]{0ex}{1.5ex}%strut
   }}
  \sum_{b=0}^{r-1} (-1)^{\delta_{i<m}} (\overline{\xi}_m)^b (\overline{\xi}_i)
   ^{-b} (\overline{m, i}) v_k 
$$
$$
\hspace{.2cm}+\frac{1}{4} \sum_{\substack{1\leq i\leq n \\ i\neq m
   \rule[0ex]{0ex}{1.5ex}%strut
   }}
   \sum_{\substack{1\leq j\leq n \\ j\neq k
   \rule[0ex]{0ex}{1.5ex}%strut
   }}
   \sum_{a,b=0}^{r-1}
  (-1)^{\delta_{i<m}+\delta_{j<k}} (\overline{\xi}_m)^b(\overline{\xi}_i)^{-b}
  (\overline{m, i})(\overline{\xi}_k)^{a}(\overline{\xi}_j)^{-a}
  (\overline{k, j}) 
$$
$$
\hspace{.4cm}
-\frac{1}{2}\sum_{\substack{1\leq i\leq n \\ i\neq m
   \rule[0ex]{0ex}{1.5ex}%strut
   }}
   \sum_{b=0}^{r-1}
   (-1)^{\delta_{i<m}} v_k(\overline{\xi}_m)^b(\overline{\xi}_i)^{-b}
  (\overline{m, i}) -\frac{1}{2} \sum_{\substack{1\leq j\leq n \\ j\neq k
   \rule[0ex]{0ex}{1.5ex}%strut
   }}
  \sum_{a=0}^{r-1} (-1)^{\delta_{j<k}}(\overline{\xi}_k)^{a}
  (\overline{\xi}_j)^{-a} (\overline{k, j}) v_m
$$
$$
\hspace{.2cm}
-\frac{1}{4} \sum_{\substack{1\leq i\leq n \\ i\neq m
   \rule[0ex]{0ex}{1.5ex}%strut
  }}\sum_{\substack
  {1\leq j\leq n \\ j\neq k
  \rule[0ex]{0ex}{1.5ex}%strut
  }}\sum_{a,b=0}^{r-1} (-1)^{\delta_{i<m}+
   \delta_{j<k}} (\overline{\xi}_k)^{a}(\overline{\xi}_j)^{-a}
   (\overline{k, j}) (\overline{\xi}_m)^b (\overline{\xi}_i)^{-b}
   (\overline{m, i})
$$

$$
%\hspace{-2.4cm}
\hspace{-.3cm}
= \frac{1}{4} \sum_{\substack{1\leq i\leq n \\ i\neq m
  \rule[0ex]{0ex}{1.5ex}%strut
  }}\sum_{\substack
  {1\leq j\leq n \\ j\neq k
  \rule[0ex]{0ex}{1.5ex}%strut
  }} \sum_{a,b=0}^{r-1} (-1)^{\delta_{i<m}
  +\delta_{j<k}} (\overline{\xi}_m)^b (\overline{\xi}_i)^{-b}
  (\overline{m, i}) (\overline{\xi}_k)^{a} (\overline{\xi}_j)^{-a}
  (\overline{k, j})
$$
$$
\hspace{.2cm}
-\frac{1}{4}\sum_{\substack{1\leq i\leq n \\ i\neq m
  \rule[0ex]{0ex}{1.5ex}%strut
   }}\sum_{\substack
   {1\leq j\leq n \\ j\neq k}}\sum_{a,b=o}^{r-1} (-1)^{\delta_{i<m}
  +\delta_{j<k}} (\overline{\xi}_k)^{a}(\overline{\xi}_j)^{-a}
  (\overline{k, j}) (\overline{\xi}_m)^b (\overline{\xi}_i)^{-b}
  (\overline{m, i})
$$
$$
\hspace{.2cm}
+\frac{1}{2} \sum_{\substack{1\leq i\leq n \\ i\neq m,k
  \rule[0ex]{0ex}{1.5ex}%strut
  }}\sum_{a,b=0}
  ^{r-1} (\overline{\xi}_m)^{a}(\overline{\xi}_k)^b (\overline{\xi}_i)^{
  -a-b} ((\overline{m, k, i})-(\overline{m, i, k})).
$$
Now in the first two of the three summands above, 
we may cancel the terms for which $\{m,i\}\cap\{k,j\}$ 
is empty. Rewriting the remaining terms and combining with the last
summand, we have the desired relation (\ref{S:bracket}).

Finally, if the $\widetilde{v}_k$ are taken instead of the
$v_k$ as generators of $H^*_{r,1,n}$ together with $G(r,1,n)$, relations
(\ref{S:trivial1})--(\ref{S:bracket}) are equivalent to the
relations for $H^*_{r,1,n}$ given in Definition \ref{S:hstar}.
\end{proof}

\begin{proof}[Proof of Theorem \ref{S:iso}]
Define an algebra homomorphism from $A_{r,1,n}$ to $H^*_{r,1,n}$ 
by sending $\overline{g}$ to
$\overline{g}$ ($g\in G(r,1,n)$) and $v_k$ to 
$\displaystyle{\frac{2}{\sqrt{3}} \ \widetilde{v}_k}$
($1\leq k\leq n$).
This map is well-defined: The relations 
(\ref{S:trivial1})--(\ref{S:action}) correspond to relations in $A_{r,1,n}$,
and (\ref{S:bracket}) corresponds to (\ref{S:89}) in $A_{r,1,n}$.
This accounts for all the relations in $A_{r,1,n}$.
The map is surjective as every generator in $H^*_{r,1,n}$ lies in its image.
On the other hand, we can define an inverse map similarly.
Therefore $A_{r,1,n}\cong H^*_{r,1,n}$.
\end{proof}

%%%%%%%%%%%%%%%%%%%%%%%%%%%%%%%%%%%%%%%%%%%%%%%%%%%%%%%%%%%%%%%%%%%%%%%
%%%%%%%%%%%%%%%%%%%%%%%%%%%%%%%%%%%%%%%%%%%%%%%%%%%%%%%%%%%%%%%%%%%%%%%
%%%%%%%%%%%%%%%%%%%%%%%%%%%%%%%%%%%%%%%%%%%%%%%%%%%%%%%%%%%%%%%%%%%%%%%

%%%%%%%%%%%%%%%%%%%%%%%%%%%%%%%%%%%%%%%%%%%%%%%%%%%%%%%%%%%%%%%555
%%%%%%%%%%%%%%%%%%%%%%%%%%%%%%%%%%%%%%%%%%%%%%%%%%%%%%%%%%%%%%%%%%%%%5
%%%%%%%%%%%%%%%%%%%%%%%%%%%%%%%%%%%%%%%%%%%%%%%%%%%%%%%%%%%%%%%%%%%55
%%%%%%%%%%%%%%%%%%%%%%%%%%%%%%%%%%%%%%%%%%%%%%%%%%%%%%%%%%%%%%%%%%%%%
%%%%%%%%%%%%%%%%%%%%%%%%%%%%%%%%%%%%%%%%%%%%%%%%%%%%%%%%%%%%%%%%%%%%%%%%5

\end{document}